\input amstex
\documentstyle{amsppt}
\magnification 1200
\vcorrection{-1cm}
\NoBlackBoxes
\input epsf

%
\def\refArtin   {1}
\def\refBLM     {2}
\def\refBKL     {3} \def\refBKLone {\refBKL}
\def\refCP      {4}
\def\refDG      {5}
\def\refEil     {6}
\def\refFH      {7}
\def\refGebGM   {8}  
\def\refGM      {9}
\def\refGMW     {10}
\def\refGL      {11}
\def\refI       {12}
\def\refKer     {13}
\def\refLM      {14}
\def\refLin     {15}
\def\refLinTalk {16}
\def\refMMP     {17}
\def\refMKS     {18}
\def\refM       {19}
\def\refUR      {20}
\def\refQPtwo   {21}
\def\refQPii    {22}


\def\sectGroup   {1.1}
\def\sectPure    {1.2}
\def\sectCable   {1.3}

\def\sectRepr    {2}

\def\sectZ      {3}
\def\sectNT     {3.1}
\def\sectCRS    {3.2}
\def\sectPerRed {3.3}

\def\sectProof   {4}
\def\sectI       {\sectProof.2}

\def\sectBfour   {5}

\def\sectGarside {6}


\def\lemAbKer   {1.1}
\def\remGL      {1.2}
\def\lemJ       {1.3}
\def\propRkJ    {1.4}

\def\lemRepr    {2.1}
\def\lemReprSix {2.2}

\def\propInvCRS {3.1}
\def\propConjCRS{3.2}
\def\propCRS    {3.3}
\def\propPerPur {3.4}
\def\propZ      {3.5}
\def\lemCable   {3.6}
\def\lemDegPer  {3.7}
\def\lemZgen    {3.8}
\def\lemNonPure {3.9}

\def\lemA     {4.1}
\def\lemInvar {4.2}
\def\lemB     {4.3}
\def\lemC     {4.4}
\def\lemD     {4.5}
\def\lemE     {4.6}
\def\lemF     {4.7}
\def\lemG     {4.8}
\def\lemConjI {4.9}
\def\lemH     {4.10}
\def\lemUniq  {4.11}
\def\lemIone  {4.12}
\def\lemItwo  {4.13}
\def\lemL     {4.14}
\def\remProof {4.15}

\def\lemZd   {\sectBfour.1}
\def\lemZall {\sectBfour.2}
\def\lemInvd {\sectBfour.3}
\def\lemInvc {\sectBfour.4}
\def\lemInvK {\sectBfour.5}
\def\lemTU   {\sectBfour.6}
\def\lemK    {\sectBfour.7}

\def\lemGarsideXY {\sectGarside.1}
\def\lemGarside   {\sectGarside.2}
\def\thQP         {\sectGarside.3}
\def\corQP        {\sectGarside.4}
\def\lemGarsideXYX{\sectGarside.5}
\def\remGarside   {\sectGarside.6}

\def\eqCommDia {1}

\def\eqLKgamma {3}
\def\eqReprA   {4}
\def\eqReprB   {5}
\def\eqLK      {6}
\def\eqLemG    {7}
\def\eqLemHa   {8}
\def\eqLemHc   {9}
\def\eqA       {10}
\def\eqB       {11}
\def\eqC       {12}
\def\eqD       {13}
\def\eqE       {14}
\def\eqF       {15}
\def\eqGL      {16}
\def\eqUTK     {17}


\def\figGraph   {1}
\def\figMove    {2}
\def\figDiamond {3}

\topmatter
\title Automorphism group of the commutator subgroup of the braid group
\endtitle
\author  S.~Yu.~Orevkov
\endauthor
\address IMT, Univ.~Paul Sabatier, Toulouse, France; Steklov Math.~Inst., Moscow, Russia    \endaddress
\endtopmatter

\def\B{\bold B}
\def\P{\bold P}
\def\J{\bold J}
\def\A{\bold A}
\def\s{\bold S}

\def\Z{\Bbb Z}
\def\C{\Bbb C}

\def\eps{\varepsilon}
\def\Aut{\operatorname{Aut}}
\def\Out{\operatorname{Out}}
\def\Sym{\operatorname{Sym}}
\def\rk {\operatorname{rk}}
\def\lk {\operatorname{lk}}
\def\id {\operatorname{id}}
\def\im {\operatorname{im}}
\def\SL {\operatorname{SL}}
\def\PSL{\operatorname{PSL}}
\def\GL {\operatorname{GL}}
\def\SC {\operatorname{SC}}
\def\card{\operatorname{Card}} \def\Card{\card}
\def\ab {{\frak{ab}}}
\def\op {{\frak{op}}}
\def\op {{\Lambda}}
\def\ad{\tilde}

\chardef\sharp=`\#
\chardef\tthat=`\^

\rightheadtext{Automorphism group of $B'_n$}
\document

\head Introduction
\endhead

Let $\B_n$ be the braid group with $n$ strings.
It is generated by $\sigma_1,\dots,\sigma_{n-1}$ (called {\it standard}
or {\it Artin} generators) subject to the relations
$$
   \text{$\sigma_i\sigma_j=\sigma_j\sigma_i$ for $|i-j|>1$;}\qquad
   \text{$\sigma_i\sigma_j\sigma_i=\sigma_j\sigma_i\sigma_j$ for
             $|i-j|=1$}.
$$
Let $\B'_n$ be the commutator subgroup of $\B_n$.
Vladimir Lin [\refLinTalk] posed a problem to compute the group of automorphisms of $\B'_n$.
In this paper we solve this problem.

\proclaim{ Theorem 1 }
If $n\ge4$, then the restriction mapping $\Aut(\B_n)\to\Aut(\B'_n)$
is an isomorphism.
\endproclaim

Dyer and Grossman [\refDG] proved that $\Out(\B_n)\cong\Z_2$ for any $n$. The only nontrivial element
of $\Out(\B_n)$ corresponds to the automorphism $\op$ defined by $\sigma_i\mapsto\sigma_i^{-1}$ for all $i=1,\dots,n$.
The center of $\B_n$ is generated by $\Delta^2$
where $\Delta=\Delta_n=\prod_{i=1}^{n-1}\prod_{j=1}^{n-i}\sigma_j$ is Garside's half-twist.
Thus $\Aut(\B_n)\cong (\B_n/\langle\Delta^2\rangle)\rtimes\Z_2$.
For an element $g$ of a group, we denote the inner automorphism $x\mapsto gxg^{-1}$ by $\ad g$.

\proclaim{ Corollary 1 } If $n\ge4$, then $\Out(\B'_n)$ is isomorphic to the dihedral group
$\bold D_{n(n-1)}=\Z_{n(n-1)}\rtimes\Z_2$.
It is generated by $\op$ and $\ad\sigma_1$
subject to the defining relations
$\op^2=\ad\sigma_1^{n(n-1)}=\op\ad\sigma_1\op\ad\sigma_1=\id$.
\endproclaim

For $n=3$, the situation is different.
It is proven in [\refGL] that $\B'_3$ is a free group of rank two generated by 
$u=\sigma_2\sigma_1^{-1}$ and $t=\sigma_1^{-1}\sigma_2$
(in fact, the free base of $\B'_3$ considered in [\refGL] is $u,v$  with $v=t^{-1}u$).
So, its automorphism group is well-known (see [\refMKS; \S3.5, Theorem N4]).
In particular (see [\refMKS; Corollary N4]),
there is an exact sequence
$$
    1\longrightarrow\B'_3\overset\iota\to\longrightarrow\Aut(\B'_3)
     \overset\pi\to\longrightarrow\GL(2,\Z)\longrightarrow 1
$$
where $\iota(x)=\tilde x$ and $\pi$ takes each automorphism of $\B'_3$ to the induced automorphism
of the abelianization of $\B'_3$ (which we identify with $\Z^2$ by choosing the images of $u$ and $t$
as a base). We have
$\ad\sigma_1(u)=t^{-1}u$, $\ad\sigma_2(u)=ut^{-1},$
$\ad\sigma_1(t)=\ad\sigma_2(t)=u$,
$\Lambda(u)=t^{-1},\;\Lambda(t)=u^{-1}$ whence
$$
    \pi(\ad\sigma_1)=\pi(\ad\sigma_2)=\left(\matrix 1 & 1 \\ -1 & 0 \endmatrix\right),
    \qquad\pi(\Lambda)=\left(\matrix 0 & -1 \\ -1 & 0 \endmatrix\right).
$$
Thus, again (as in the case $n\ge 4$) the image of $\Aut(\B_3)$ in  $\Out(\B'_3)\cong\GL(2,\Z)$ is isomorphic
to $\Bbb D_6$ but this time it is not the whole group $\Out(\B'_3)$.

Let $\s_n$ be the symmetric group and $\A_n$ its alternating subgroup.
Let $\mu=\mu_n:\B_n\to\s_n$ be the homomorphism which takes $\sigma_i$ to the
transposition $(i,i+1)$ and let $\mu'$ be the restriction of $\mu$ to $\B'_n$.
Then $\P_n=\ker\mu_n$ is the group of {\it pure braids}.
Let $\J_n = \B'_n\cap \P_n = \ker\mu'$. Note that the image of $\mu'$ is $\A_n$.
The following diagram commutes where the 
rows are exact sequences and all the
unlabeled
arrows (except ``$\to1$")
are inclusions:
$$
\matrix
1 & \longrightarrow &\J_n &\longrightarrow &\B'_n &\overset{\mu'}\to\longrightarrow& \A_n&\longrightarrow &1\\
  &             &\overset{|}\to\downarrow\;\,& &\overset{|}\to\downarrow\;\,& &\overset{|}\to\downarrow\;\,&\\
1 & \longrightarrow &\P_n &\longrightarrow &\B_n &\overset{\mu}\to\longrightarrow& \s_n&\longrightarrow &1
\endmatrix
                                                                     \eqno(\eqCommDia)
$$

Recall that a subgroup of a group $G$ is called {\it characteristic} if it is invariant
under each automorphism of $G$.
Lin proved in [\refLin; Theorem D] that $\J_n$ is a characteristic subgroup of $\B'_n$ for $n\ge 5$
(note this fact is used in our proof of Theorem 1 for $n\ge 5$).
By Theorem 1, this result extends to the case $n=4$.

\proclaim{ Corollary 2 } $\J_4$ is a characteristic subgroup of $\B'_4$.
\endproclaim

Note that $\J_3$ is not a characteristic subgroup of $\B'_3$. Indeed, let $\varphi\in\Aut(\B'_3)$
be defined by $u\mapsto u$, $t\mapsto ut$. Then $ut\in\J_3$ whereas $\varphi^{-1}(ut)=t\not\in\J_3$.


\head 1. Preliminaries
\endhead

Let $e:\B_n\to\Z$ be the homomorphism defined by $e(\sigma_i)=1$ for all $i=1,\dots,n$.
Then we have $\B'_n = \ker e$.

\subhead\sectGroup. Groups
\endsubhead
For a group $G$, we denote its unit element by $1$, the center by $Z(G)$,
the commutator subgroup by $G'$, the second commutator subgroup $(G')'$ by $G''$, and
the abelianization $G/G'$ by $G^\ab$. We denote $x^{-1}yx$ by $y^x$
(thus $\tilde x(y^x)=y$) and we denote
the commutator $xyx^{-1}y^{-1}$ by $[x,y]$.
For $g\in G$, we denote the centralizer of $g$ in $G$ by $Z(g,G)$.
If $H$ is a subgroup of $G$, then, evidently, $Z(g,H)=Z(g,G)\cap H$.

\proclaim{ Lemma \lemAbKer }
Let $G$ be a group generated by a set $\Cal A$. Assume that there exists a
homomorphism $e:G\to\Bbb Z$ such that $e(\Cal A)=\{1\}$.
Let $\bar e$ be the induced homomorphism $G^\ab\to\Z$.
Let $\Gamma$ be the graph such that the set of vertices is $\Cal A$ and two vertices
$a$ and $b$ are connected by an edge when $[a,b]=1$.

If the graph $\Gamma$ is connected,
then $(\ker e)^\ab \cong \ker\bar e$.
\endproclaim

\demo{ Proof }
Let $K=\ker e$. 
Let us show that $G'\subset K'$.
Since $K'$ is normal in $G$, and $G'$ is the normal closure  of
the subgroup generated by
$[a,b]$, $a,b\in\Cal A$, it is enough to show that $[a,b]\in K'$ for any $a,b\in\Cal A$.

We define a relation $\sim$ on $\Cal A$ by setting
$a\sim b$ if $[a,b]\in K'$. Since $\Gamma$ is connected, it remains to note that this
relation is transitive. Indeed, if $a\sim b\sim c$, then
$[a,c]=[a,b][bab^{-2},bcb^{-2}][b,c]\in K'$.

Thus $G'\subset K'$ whence $G'=K'$ and we obtain $K^\ab=K/K'=K/G'=\ker\bar e$.
%
\qed\enddemo

\medskip\noindent{\bf Remark \remGL.}
The fact that $\B_n''=\B_n'$ for $n\ge 5$ proven by Gorin and Lin [\refGL] (see also [\refLin; Remark 1.10])
is an immediate corollary of Lemma \lemAbKer. Indeed, if we set $G=\B_n$ and
$\Cal A=\{\sigma_i\}_{i=1}^{n-1}$, then $\Gamma$ is connected whence
$\B_n'/\B''_n =(\ker e)^\ab \cong \ker\bar e = \{1\}$.
In the same way we obtain $G''=G'$ when $G$ is Artin group of type $D_n$ ($n\ge 5$),
$E_6$, $E_7$, $E_8$, $F_4$, or $H_4$.
\medskip

\subhead\sectPure. Pure braids
\endsubhead
Recall that $\P_n$ is generated by the braids $\sigma_{ij}^2$, $1\le i<j\le n$, where
$\sigma_{ij}=\sigma_{ji}=\sigma_{j-1}\dots\sigma_{i+1}\sigma_i\sigma_{i+1}^{-1}\dots\sigma_{j-1}^{-1}$.
For a pure braid $X$, let us denote the linking number of the $i$-th and $j$-th strings by
$\lk_{ij}(X)$. If $X$ is presented by a diagram with under- and over-crossings, then
$\lk_{ij}(X)$ is the half-sum of the signs of those crossings where the $i$-th and $j$-th strings cross. 
Let $A_{ij}$ be the image of $\sigma_{ij}^2$ in $\P_n^\ab$. We have, evidently,
$$
         \lk_{\gamma(i),\gamma(j)}(X) = \lk_{i,j}(X^\gamma),
             \qquad\text{for any $X\in\P_n$, $\gamma\in\B_n$}                \eqno(\eqLKgamma)
$$
(here $\gamma(i)=\mu(\gamma)(i)$ which is coherent with the interpretation of $\B_n$ with a
mapping class group; see \S\sectNT).

It is well known that $\P_n^\ab$ is freely generated by $\{A_{ij}\}_{1\le i<j\le n}$.
This fact is usually derived from Artin's
presentation of $\P_n$ (see [\refArtin; Theorem 18]) but it also admits a very simple self-contained proof
based on the linking numbers. Namely, let $L$ be free abelian group with a free base
$\{a_{ij}\}_{1\le i<j\le n}$. Then it is immediate to check that
the mapping $\P_n\to L$, $X\mapsto\sum_{i<j}\lk_{i,j}(X)a_{ij}$ is a homomorphism and that the
induced homomorphism $\P_n^\ab\to L$ is the inverse of $L\to\P_n^\ab$, $a_{ij}\mapsto A_{ij}$.
In particular, we see that the quotient map $\P_n\to\P_n^\ab$ is given by
$X\mapsto\sum_{i<j}\lk_{i,j}(X)A_{ij}$.

\proclaim{ Lemma \lemJ } If $n\ge 5$, then the mapping $\J_n\to\P_n^\ab$,
$X\mapsto\sum\lk_{ij}(X)A_{ij}$ defines an isomorphism
$\J_n^\ab\cong\{\sum x_{ij}A_{ij}\mid\sum x_{ij}=0\}\subset\P_n^\ab$. \qed
\endproclaim

\demo{ Proof }
Follows from Lemma \lemAbKer\ with $\P_n$, $e|_{\P_n}$, and $\{\sigma_{ij}^2\}_{1\le i<j\le n}$
standing for $G$, $e$, and $\Cal A$ respectively.
\qed
\enddemo

So, when $n\ge 5$, we identify $\J_n^\ab$ with its image in $\P_n^\ab$.
The following fact will not be used in the proof of Theorems 1.

\proclaim{ Proposition \propRkJ }
(a). $\J_n^\ab$ is a free abelian group and
$$
     \rk \J_n^\ab = \left(\matrix n\\2\endmatrix\right) + \cases  \;\;\,1, & n\in\{3,4\},\\
                        -1, & \text{otherwise}.
     \endcases
$$
(b). 
$\Cal E_3=\{\bar u\bar t,\bar t\bar u,\bar u^3,\bar t^3\}$
and $\Cal E_4=\Cal E_3\cup\{\bar c^2,\bar w^2,(\bar c\bar w)^2\}$
are free bases of $\J_3^\ab$ and $\J_4^\ab$ respectively.

\smallskip

{\rm(Here and below $\bar x$ stands for the image of $x$ under the quotient map $\J_n\to\J_n^\ab$.)}

\smallskip\noindent
(c). Let $p_n:\J_n^\ab\to\P_n^\ab$, $n=3,4$, be induced by $\J_n\to\P_n\to\P_n^\ab$.
Then
$$
    \im p_n = \{\sum x_{ij}A_{ij}\mid\sum x_{ij}=0\},
    \qquad \ker p_n = \langle\bar u^3, \bar t^3 \rangle.
$$
\endproclaim

\demo{ Proof }
(a). The result is obvious for $n=2$ and it follows from Lemma \lemJ\ for $n\ge 5$.

For $n=3$, the result follows from the following argument proposed by the referee.
We have $\B'_3\cong\pi_1(\Gamma)$ where $\Gamma$ is the bouquet $S^1\vee S^1$. Since $|\B'_3/\J_3|=3$ 
(see (\eqCommDia)), we have $\J_3^\ab\cong H_1(\tilde\Gamma)$ where $\tilde\Gamma\to\Gamma$
is a connected $3$-fold covering. Then the Euler characteristic of $\tilde\Gamma$ is
$\chi(\tilde\Gamma)=3\chi(\Gamma)=-3$ whence $\rk H_1(\tilde\Gamma)=4$.


The group $\J_4^\ab$ can be easily computed by Reidemeister--Schreier method
either as $\ker\mu'$ using Gorin and Lin's [\refGL] presentation for $\B'_4$, or
as $\ker(e|_{\P_4})$ using Artin's presentation [\refArtin] of $\P_4$.
Here is the GAP code for the first method:

\smallskip
\vbox
{\tt\noindent
f:=FreeGroup(4); u:=f.1; v:=f.2; w:=f.3; c:=f.4;

\noindent
g:=f/[u*c/u/w, u*w/u/w*c/w/w, v*c/v/w*c, v*w/v/w*c*c/w*c/w*c/w*c]; 

\noindent
u:=g.1; v:=g.2; w:=g.3; c:=g.4; \ \ \ \sharp\ group B'(4) according to [\refGL]

\noindent
s:=SymmetricGroup(4); t1:=(1,2); t2:=(2,3); t3:=(3,4);

\noindent
U:=t2*t1; V:=t1*t2; W:=t2*t3*t1*t2; C:=t3*t1; \sharp\ U=mu(u),V=mu(v),...

\noindent
mu:=GroupHomomorphismByImages(g,s,[u,v,w,c],[U,V,W,C]);

\noindent
AbelianInvariants(Kernel(mu)); \ \ \ \ \sharp\ should be [0,0,0,0,0,0,0]
}

\midinsert
\centerline{\epsfxsize=55mm\epsfbox{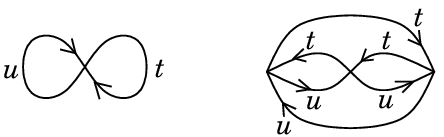}}
\botcaption{ Figure \figGraph } The graphs $\Gamma$ and $\tilde\Gamma$ in the proof of Proposition \propRkJ
\endcaption
\endinsert

\smallskip
(b) for $n=3$. In Figure \figGraph\ we show the graphs $\Gamma$ and $\tilde\Gamma$ discussed above.
We see that the loops in $\tilde\Gamma$ represented
by the elements of $\Cal E_3$ form a base of $H_1(\tilde\Gamma)$.

\smallskip
(c) for $n=3$. The claim about $\im p_3$ is evident and a computation of the linking numbers shows that
$p_3(\bar u^3)=p_3(\bar t^3)=0$.
Note that the braid closures of both $u^3$ and $t^3$ are Borromean links.
So, maybe, it could be interesting to study how all this stuff is related to Milnor's $\mu$-invariant.


\smallskip
(b,c) for $n=4$. The claim about $\im p_4$ is evident and a computation of the linking numbers shows
that $p_4(\Cal E_4\setminus\{\bar t^3,\bar u^3\})$ is a base of $\im p_4$.
One can check that the homomorphism
$\B'_4\to\B'_4/\bold K_4\cong\B'_3$ maps $\J_4$ to $\J_3$. Hence it induces a homomorphism
$\J_4^\ab\to\J_3^\ab$ which takes $\bar u^3$ and $\bar t^3$ of $\J_4^\ab$ to
$\bar u^3$ and $\bar t^3$ of $\J_3^\ab$. Hence 
$\rk(\ker p_4)\ge\rk\langle\bar u^3,\bar t^3\rangle=2$.
Since $\rk\J_4^\ab=7$ and $\rk(\im p_4)=5$, we conclude that 
$\ker p_4=\langle\bar u^3,\bar t^3\rangle$.
%
\qed
\enddemo



\subhead\sectCable. Mixed braid groups and the cabling map
\endsubhead
Let $n\ge 1$ and $\vec m=(m_1,\dots,m_k)$, $m_1+\dots+m_k=n$, $m_i\in\Z$, $m_i>0$.

The {\it mixed braid group} $\B_{\vec m}$ (see [\refM], [\refUR], [\refGMW])
is defined as $\mu^{-1}(S_{\vec m})$ where
$S_{\vec m}$ is the stabilizer of the following vector under the natural action of $\s_n$ on $\Z^n$:
$$
    (\,\underset{m_1}\to{\underbrace{1,\dots,1}}\,,\underset{m_2}\to{\underbrace{2,\dots,2}},\,
     \dots,\,\underset{m_k}\to{\underbrace{k,\dots,k}}\,).
$$
We emphasize two particular cases: $\B_{1,\dots,1}$ is the pure braid group and
$\B_{n-1,1}$ is the Artin group corresponding to the Coxeter group of type $B_{n-1}$.

We define the
{\it cabling map}
$\psi=\psi_{\vec m}:\B_k\times(\B_{m_1}\times\dots\times\B_{m_k})\to\B_n$
by sending $(X;X_1,\dots,X_k)$ to the braid obtained by replacing each strand of $X$
by a geometric braid representing $X_i$ embedded into a small tubular neighbourhood
of this strand. 

Note that $\psi_{\vec m}$ is not a homomorphism but its restriction to
$\P_k\times\prod_i\B_{m_i}$ is.
We have $\psi(\P_k\times\prod_i\P_{m_i})\subset\P_n$ and
$\psi(\P_k\times\prod_i\B_{m_i})\subset\B_{\vec m}$.



\head\sectRepr. $\J_n^\ab$ as an $\A_n$-module and its automorphisms
\endhead
Let $n\ge 5$.
As we mentioned already, by [\refLin; Theorem D], $\J_n$ is a characteristic subgroup of $\B'_n$, i.e., $\J_n$
is invariant under any automorphism of $\B'_n$ (in fact a stronger statement is proven in [\refLin]).

\proclaim{ Lemma \lemRepr }
Let $\varphi\in\Aut(\B'_n)$ be such that $\mu'\varphi = \mu'$.
Let $\varphi_*$ be the automorphism of $\J_n^\ab$ induced by $\varphi|_{\J_n}$.
Then $\varphi_* = \pm\id$.
\endproclaim

\demo{ Proof }
The exact sequence $1\to\J_n\to\B'_n\to\A_n\to 1$ (see (\eqCommDia)) defines
an action of $\A_n$ on $\J_n^\ab$ by conjugation.
The condition $\mu'\varphi=\mu'$ implies that $\varphi_*$ is $\A_n$-equivariant.
Let $V$ be a complex vector space with base $e_1,\dots,e_n$ endowed with
the natural action of $\s_n$ induced by the action on the base.
We identify $\P_n^\ab$ with its image
in the symmetric square $\Sym^2 V$ by the homomorphism $A_{ij}\mapsto e_i e_j$.
Then, by Lemma \lemJ, we may identify $\J_n^\ab$ with $\{\sum c_{ij}e_ie_j\mid\sum c_{ij}=0\}$.
These identifications are compatible with the action of $\A_n$.

For a partition $\lambda=(\lambda_1,\dots,\lambda_r)$, we denote the corresponding
irreducible representation of $\s_n$ over $\C$ (the $\C\s_n$-module) by $V_\lambda$, see, e.g., [\refFH; \S4].
For an element $v$ of a $\Bbb C\s_n$-module, let $\langle v\rangle_{\C\s_n}$ be the
$\C\s_n$-submodule generated by $v$.
We set $e_0=e_1+\dots+e_n$, $U=\langle e_0\rangle_{\C\s_n}=\Bbb C e_0$, and $U^\perp=\langle e_1-e_2\rangle_{\C\s_n}$.
Consider the following $\Bbb C\s_n$-submodules of $\Sym^2 V$:
$$
\split
   W_0&=\langle e_1^2\rangle_{\C\s_n},\quad
   W_1 =\langle w\rangle_{\C\s_n}=\C w\quad\text{where $w=\sum e_i e_j$},\\
   W_2&=\langle(e_1-e_2)(e_3+\dots+e_n)\rangle_{\C\s_n},\quad
   W_3=\langle(e_1-e_2)(e_3-e_4)\rangle_{\C\s_n}.
\endsplit
$$
We have $\Sym^2 V=\Sym^2(U\oplus U^\perp)=\Sym^2 U\oplus\Sym^2 U^\perp \oplus(U\otimes U^\perp)$
and $U^\perp\cong V_{n-1,1}$ (that is $V_\lambda$ for $\lambda=(n-1,1)$).
It is known (see [\refMMP; Lemma 2.1] or [\refFH; Exercise 4.19]) that
$\Sym^2 V_{n-1,1}\cong U\oplus V_{n-1,1}\oplus V_{n-2,2}\cong V\oplus V_{n-2,2}$. Thus
$$
   \Sym^2 V \cong  V\oplus V\oplus V_{n-2,2}.                                    \eqno(\eqReprA)
$$
Let $W=\J_n^\ab\otimes\C$. It is clear that $\Sym^2 V=W_0\oplus W_1\oplus W$.
Since $W_0\cong V$ and $W_1\cong U$, we obtain $W\cong U^\perp\oplus V_{n-2,2}$ by cancelling out
$U\oplus V$ in (\eqReprA).
Note that $(e_1-e_2)(e_3+\dots+e_n)=(e_1-e_2)(e_0-(e_1+e_2))=(e_1-e_2)e_0-(e_1^2-e_2^2)$,
hence the mapping  $e_i-e_j\mapsto (e_i-e_j)e_0-(e_i^2-e_j^2)$ induces
an isomorphism of $\C\s_n$-modules $U^\perp\cong W_2$. The identity
$$
   (n-2)(e_1-e_2)e_3 = (e_1-e_2)(e_3+\dots+e_n) +  \sum_{i\ge 4} (e_1-e_2)(e_3-e_i) \eqno(\eqReprB)
$$
shows that $W_2+W_3 = \langle(e_1-e_2)e_3\rangle_{\C\s_n}=W$. One easily checks that
$W_2$ and $W_3$ are orthogonal to each other with respect to the scalar product on $W+W_1$ for which
$\{e_ie_j\}_{i,j}$ is an orthonormal basis. Therefore $W=W_2\oplus W_3$ is the decomposition of
$W$ into irreducible factors.

We have $W_2\cong V_{n-1,1}$ and $W_3\cong V_{n-2,2}$. Since the corresponding Young diagrams are
not symmetric, $W_2$ and $W_3$ are irreducible as $\C\A_n$-modules (see [\refFH; \S 5.1]).
Since $\dim W_2\ne\dim W_3$ 
and $\varphi_*$ is $\A_n$-equivariant,
Schur's lemma implies that $\varphi_*|_{W_k}$, $k=2,3$, is
multiplication by a constant $c_k$. Moreover, since $\varphi_*$ is an automorphism of $\J_n^\ab$
(a discrete subgroup), we have $c_k=\pm1$. If $c_3=-c_2=\pm1$, then (\eqReprB) contradicts
the fact that $\varphi_*((e_1-e_2)e_3)\in\J_n^\ab$.
\qed\enddemo


Let $\nu\in\Aut(\s_6)$ be defined by $(12)\mapsto(12)(34)(56)$, $(123456)\mapsto(123)(45)$.
It is well known that $\nu$ represents the only nontrivial element of $\Out(\s_6)$.

\proclaim{ Lemma \lemReprSix } Let $\varphi\in\Aut(\B'_6)$.
Then $\mu'\varphi \ne \nu\mu'$.
\endproclaim

\demo{ Proof }
Given a commutative ring $k$ and a $k\A_6$-module
$V$ corresponding to a representation $\rho:\A_6\to\GL(V,k)$, we denote the $k\A_6$-module
corresponding to the representation $\rho\nu$ by $\nu^*(V)$.
It is clear that $\nu^*$ is a covariant functor which preserves direct sums
(hence irreducibility), tensor products, symmetric powers etc.

Suppose that $\mu'\varphi = \nu\mu'$.
As in the proof of Lemma \lemRepr, we endow $\J_6^\ab$ with the action of $\A_6$.
The condition $\mu'\varphi=\nu\mu'$ implies that $\varphi$ induces an isomorphism
of $\A_6$-modules $\J_6^\ab\cong \nu^*(\J_6^\ab)$.
Let us show that these modules are not isomorphic.

We have $\J_6^\ab\otimes\C\cong V_{5,1}\oplus V_{4,2}$ (see the proof of Lemma \lemRepr).
Hence $\nu^*(\J_6^\ab)\otimes\C\cong\nu^*(V_{5,1})\oplus\nu^*(V_{4,2})$.
We have $\dim V_{5,1}=5\ne 9=\dim V_{4,2}$, thus, to complete the proof, it is enough to show that
$V_{5,1}\not\cong\nu^*(V_{5,1})$ (note that $V_{4,2}\cong\nu^*(V_{4,2})$).
Indeed, $\nu$ exchanges the conjugacy classes of the permutations
$a=(123)$ and $b=(123)(456)$, hence we have $\chi(a)=2\ne-1=\chi(b)=\chi\nu(a)$ where
$\chi$ and $\chi\nu$ are the characters of $\A_6$ corresponding to $V_{5,1}$ and to
$\nu^*(V_{5,1})$ respectively.
\qed\enddemo



\head\sectZ. Centralizers of pure braids
\endhead

Centralizers of braids are computed by Gonz\'alez-Meneses and Wiest [\refGMW].
For pure braids the answer is much simpler and it can be easily obtained as a specialization
of the results of [\refGMW].


\subhead\sectNT. Nielsen-Thurston trichotomy
\endsubhead
The following definitions and facts we reproduce from [\refGMW; Section 2] where they are
taken from different sources, mostly from
the book [\refI] which can be also used as a general introduction to the subject.

Let $\Bbb D$ be a disk in $\Bbb C$ that contains $X_n=\{1,\dots,n\}$.
The elements of $X_n$ will be called {\it punctures}.
It is well known that $\B_n$ can be identified with the mapping class group
$\Cal D/\Cal D_0$ where $\Cal D$ is the group of diffeomorphisms $\beta:\Bbb D\to\Bbb D$
such that $\beta|_{\partial\Bbb D}=\operatorname{id}_{\partial\Bbb D}$ and $\beta(X_n)=X_n$,
and $\Cal D_0$ is the connected component of the identity.
Sometimes, by abuse of notation, we shall not distinguish between braids and
elements of $\Cal D$.
For $A,B\subset\Bbb D$, we write $A\sim B$
if $\beta_0(A)=B$ for some $\beta_0\in\Cal D_0$.

An embedded circle in $\Bbb D\setminus X_n$ is called an {\it essential curve}
if it encircles more than one but less than $n$ points of $X_n$.
A {\it multicurve} in $\Bbb D\setminus X_n$ is a disjoint union of embedded circles.
It is called {\it essential} if all its components are essential.

Let $\beta\in\Cal D$.
We say that a multicurve $C$ in $\Bbb D\setminus X_n$ is {\it stabilized} or {\it preserved} by $\beta$
if $\beta(C)\sim C$
(the components of $C$ may be permuted by $\beta$).
The braid represented by $\beta$
is called {\it reducible} if $\beta$ stabilizes some essential multicurve. 

A braid $\beta$ is called {\it periodic} if some power of $\beta$ belongs to $Z(\B_n)$.
If a braid is neither periodic nor reducible, then it is called {\it pseudo-Anosov}; see [\refI].


\subhead\sectCRS. Canonical reduction systems. Tubular and interior braids
\endsubhead
An essential curve $C$ is called a {\it reduction curve} for a braid $\beta$ if
it is stabilized by some power of $\beta$ and any other curve stabilized by some power of $\beta$
 is isotopic in $\Bbb D\setminus X_n$ to a curve disjoint from $C$.
An essential multicurve is called a {\it canonical reduction system} (CRS) for $\beta$
if its components represent all isotopy classes of reduction curves for $\beta$
(each class being represented once).
%
It is known that there exists a canonical reduction system for any braid
and that it is unique up to isotopy, see [\refBLM], [\refI; \S7], [\refGMW; \S2].
If a braid is periodic or pseudo-Anosov, the CRS is empty.
The following properties of CRS are immediate consequences from the existence and uniqueness.

\proclaim{ Proposition \propInvCRS } Let $C$ be the CRS for $\beta\in\Cal D$.
Then $C$ is the CRS for $\beta^{-1}$. \qed
\endproclaim

\proclaim{ Proposition \propConjCRS }
Let $\beta,\gamma\subset\Cal D$ and let $C$ be the CRS for $\beta$.
Then $\gamma^{-1}(C)$ is the CRS for $\beta^\gamma$. \qed
\endproclaim

\proclaim{ Proposition \propCRS } 
Let $\beta,\gamma\in\Cal D$ represent
commuting braids. Then:

\smallskip
(a). $\gamma$ preserves the CRS of $\beta$.

\smallskip
(b). If $\gamma$ is pure, then it preserves each reduction curve of $\beta$.
\endproclaim

\demo{ Proof }
(a). Follows from Proposition \propConjCRS.

(b). Follows from (a).
\qed\enddemo

We say that a braid is in {\it almost regular form} if its CRS
is a union of round circles (`almost' because the definition of regular form in [\refGMW] includes
some more conditions which we no not need here).
By Proposition \propConjCRS\ any braid has a conjugate in almost regular form.

%
%
%
%
%

Let $\beta$ be an element of $\Cal D$ which represents a reducible braid in almost regular form and
let $C$ be a CRS for $\beta$. Without loss of generality we may assume that $\beta(C)=C$
and $C$ is a union of round circles.
Let $R'$ be the union of the outermost components of $C$.
We add to $R'$ small circles around the points of $X_n$ not encircled by curves from $R'$
and we denote the resulting multicurve by $R$.
Let $C_1,\dots,C_k$ be the connected components of $R$ numbered from the left to the right.

Recall that the geometric braid (a union of strings in the cylinder $[0,1]\times\Bbb D$)
is obtained from $\beta$ as follows. Let $\{\beta_t:\Bbb D\to\Bbb D\}_{t\in[0,1]}$ be an isotopy such that
$\beta_0=\beta$, $\beta_1=\id_{\Bbb D}$, and $\beta_t|_{\partial\Bbb D}=\id_{\partial\Bbb D}$ for any $t$.
Then the $i$-th string of the geometric braid is the graph of the mapping $t\mapsto\beta_t(i)$
and the whole geometric braid is $\bigcup_t\big(\{t\}\times\beta_t(X_n)\big)$.
Similarly, starting from the circles $C_i$, we define the embedded cylinders (tubes)
$\bigcup_t\big(\{t\}\times\beta_t(C_i)\big)$, $i=1,\dots,k$.

Let $m_i$ be the number of punctures encircled by $C_i$.
Following [\refGMW; \S 5.1], we define the {\it interior braid} $\beta_{[i]}\in\B_{m_i}$, $i=1,\dots,k$,
as the element of $\B_{m_i}$ corresponding to the union of strings contained in the $i$-th tube,
and we define the {\it tubular braid} $\hat\beta$ of $\beta$ as the braid obtained by shrinking each tube
to a single string. Let $\vec m=(m_1,\dots,m_k)$ and let $\psi_{\vec m}$ be the
cabling map (see \S\sectCable). Then we have
$\beta = \psi_{\vec m}(\hat\beta;\beta_{[1]},\dots,\beta_{[k]})$.

Recall that $C$ is a CRS for $\beta$.
Let $a$ be an open connected subset of $\Bbb D$ such that $\partial a\subset C\cup\partial\Bbb D$.
To each such $a$ we associate the braid which is the union of the strings of $\beta$ starting
at $a$ and the strings obtained by shrinking the tubes corresponding to the interior components of $\partial a$.
We denote this braid by $\beta_{[a]}$. 
For example, if $a$ is the exterior component of $\Bbb D\setminus C$, then $\beta_{[a]}=\hat\beta$.
%



\subhead\sectPerRed. 
Periodic and reducible pure braids
\endsubhead
The structure of the centralizers of
periodic and reducible braids becomes extremely simple if we restrict
our attention to pure braids only.
The following fact immediately follows from  a result due to
Eilenberg [\refEil] and Ker\'ekj\'art\'o [\refKer]
(see [\refGMW; Lemma 3.1]).

\proclaim{ Proposition \propPerPur }
A pure braid is periodic if and only if it is a power of $\Delta^2$.
\endproclaim

The following fact can be considered as a specialization of the results of [\refGMW].

\proclaim{ Proposition \propZ } Let $\beta$ be a pure $n$-braid.

\smallskip
(a). If $\beta$ is periodic, then $Z(\beta;\P_n)=\P_n$.

\smallskip
(b). If $\beta$ is pseudo-Anosov, then $Z(\beta;\P_n)$ is the free abelian group
generated by $\Delta^2$ and some pseudo-Anosov braid which may or may not
coincide with $\beta$.

\smallskip
(c). If $\beta$ is reducible non-periodic and in almost regular form, then
$\psi_{\vec m}$ maps isomorphically
$Z(\hat\beta;\P_k)\times Z(\beta_{[1]};\P_{m_1})\times\dots\times Z(\beta_{[k]};\P_{m_k})$
onto $Z(\beta;\P_n)$ {\rm(see  \S\sectCRS)}.
\endproclaim

\demo{ Proof }
(a). Follows from Proposition \propPerPur.

(b). Follows from [\refGMW; Proposition 4.1].

(c). (See also the proof of [\refGMW; Proposition 5.17]).
By Proposition \propCRS\ we have $Z(\beta;\P_n)\subset\psi_{\vec m}(\P_k\times\prod\P_{m_i})$.
The injectivity of the considered mapping and the fact that $\psi_{\vec m}^{-1}(Z(\beta;\P_n))$
is as stated, are immediate consequences from the following observation: if two geometric braids are isotopic,
then the braids obtained from them by removal of some strings are isotopic as well.
\qed\enddemo

\proclaim{ Lemma \lemCable } Let $\vec m=(m_1,\dots,m_k)$, $m_1+\dots+m_k=n$, and $p\in\Z$. Then
$\psi_{\vec m}(\Delta^{p}_k;\Delta_{m_1}^{p},\dots,\Delta_{m_k}^{p})=\Delta_n^{p}$.
\endproclaim

\demo{ Proof } The result immediately follows from the geometric characterization of $\Delta$ as
a braid all whose strings lie on a half-twisted band. Note that the sub-bands of the
half-twisted band arising from consecutive strings also consist of half-twisted bands.
\qed\enddemo

If $X$ is a periodic pure braid, then $X=\Delta^{2d}$, $d\in\Z$, by Proposition \propPerPur.
In this case we set $d=\deg X$, the {\it degree} of $X$. It is clear that
$\lk_{ij}(X)=d$ for any $i<j$.

\proclaim{ Lemma \lemDegPer } Let $C$ be the CRS for a reducible pure braid represented by $\beta\in\Cal D$.
Let $a$ and $b$ be two neighboring components of $\Bbb D\setminus C$ and let $X=\beta_{[a]}$ and $Y=\beta_{[b]}$
be the braids associated to $a$ and $b$ {\rm(}see the end of \S\sectCRS\/{\rm)}.
Suppose that each of $X$ and $Y$ is periodic.
Then $\deg X\ne\deg Y$.
\endproclaim

\demo{ Proof } Suppose that $\deg X = \deg Y = p$, i.~e., $X=\Delta_k^{2p}$ and $Y=\Delta_m^{2p}$
for some $k,m\ge2$.
Let $C_i$ be the component of $C$ that separates $a$ and $b$.
We may assume that $a$ is exterior to $C_i$. Let $c$ be the closure of $a\cup b$.
Then we have 
$$
   \beta_{[c]} = \psi_{1,\dots,1,m,1,\dots,1}
         (\Delta_k^{2p};1,\dots,1,\Delta_m^{2p},1,\dots,1)=\Delta_{k+m-1}^{2p}
$$
by Lemma \lemCable. Hence $\beta_{[c]}$ preserves any closed curve, in particular a curve
which separates some two strings of $\beta_{[b]}$ and encircles a string of $\beta_{[c]}$ not
belonging to $\beta_{[b]}$. Such a curve is not isotopic to any curve disjoint from $C_i$.
This fact contradicts the condition that $C_i$ is a reduction curve.
\qed\enddemo

\proclaim{ Lemma \lemZgen } $Z(\sigma_1^2\sigma_3^{-2};\J_n)\cong\P_{n-2}\times\Z$ for $n\ge 4$.
\endproclaim

\demo{ Proof } The CRS for $\sigma_1^2\sigma_3^{-2}$ consists of two round circles:
one of them encircles the punctures 1 and 2, and the other one encircles
the punctures 3 and 4.
Then Proposition \propPerPur(c) implies that
$\psi=\psi_{\vec m}:\P_{n-2}\times(\P_2\times\P_2)\to\P_n$, $\vec m=(2,2,1_{n-4})$,
is injective and $\im\psi=Z(\sigma_1^2\sigma_3^{-2};\P_n)$. 
One easily checks that 
the mapping $\P_{n-2}\times\P_2\to Z(\sigma_1^2\sigma_3^{-2};\J_n)$,
$(X,\sigma_1^{k})\mapsto\psi(X;\sigma_1^{k},\sigma_1^{-m})$, $m=e(\psi(X;1,1))+k$ is an isomorphism.
Indeed, any element of $Y\in Z(\sigma_1^2\sigma_3^{-2};\P_n)$
is of the form $Y=\psi(X;\sigma_1^{k},\sigma_1^{-m})$ and the condition $e(Y)=0$
reads as $e(\psi(X;1,1))+k-m=0$.
\qed\enddemo

\proclaim{ Lemma \lemNonPure }
Let $\beta$ be a reducible $n$-braid in almost regular form.
Suppose that $\hat\beta\in\P_k$ and that $\beta_{[i]}$'s are pairwise
non-conjugate to each other {\rm(see \S\sectCRS)}. Then $\psi_{\vec m}$ maps isomorphically
$Z(\hat\beta;\P_k)\times Z(\beta_{[1]};\B_{m_1})\times\dots\times Z(\beta_{[k]};\B_{m_k})$
onto $Z(\beta;\B_n)$
\endproclaim

\demo{ Proof } See the proof of Proposition \propZ(c).
\qed\enddemo



\head\sectProof. Proof of Theorem 1 for $n\ge5$
\endhead

\subhead\sectProof.1. Invariance of the conjugacy class of $\sigma_1\sigma_3^{-1}$
\endsubhead
Suppose that $n\ge 5$. Let $\varphi\in\Aut(\B'_n)$ be such that $\mu'\varphi=\mu'$
and $\varphi_*=\id$ where $\varphi_*$ is as in
Lemma \lemRepr. Then  
we have
$$
    \lk_{ij}(X) = \lk_{ij}(\varphi(X)),\qquad X\in\J_n, \quad 1\le i<j\le n.              \eqno(\eqLK)
$$

Let $\tau = \psi_{2,n-2}(1;\sigma_1^{(n-2)(n-3)},\Delta^{-2})$. We have $\tau\in\J_n$.

\proclaim{ Lemma \lemA } Let $X$ be $\sigma_1^{2}\sigma_3^{-2}$, $k\ne0$, or $\tau$.
Let $\alpha\in\Cal D$ represent $\varphi(X)$.
Let $C$ be a simple closed curve preserved by $\alpha$.
Suppose that $C$ encircles at least two punctures.
Then the punctures $1$ and $2$ are in the same
component of $\Bbb D\setminus C$.
\endproclaim

\demo{ Proof } Suppose that $1$ and $2$ are separated by $C$. Without loss of generality we
may assume that $1$ is outside $C$ and $2$ is inside $C$.
Let $p$ be another puncture inside $C$.
Then we have $\lk_{1,p}(\alpha)=\lk_{1,2}(\alpha)$ which contradicts (\eqLK) because
$\lk_{1,2}(X)\ne0$ and $\lk_{1,p}(X)=0$ for any $p\ne2$. \qed
\enddemo

\proclaim{ Lemma \lemInvar }
Let $\alpha\in\Cal D$ represent $\varphi(\sigma_1^2\sigma_3^{-2})$.
Then the CRS for $\alpha$ is invariant under 
some element of $\Cal D$ which exchanges $\{1,2\}$ and $\{3,4\}$.
\endproclaim

\demo{ Proof } Follows from Propositions \propInvCRS\ and \propConjCRS\ because $\alpha$
is conjugate to $\alpha^{-1}$ and the conjugating element of $\Cal D$ exchanges
 $\{1,2\}$ and $\{3,4\}$.
\qed\enddemo

\proclaim{ Lemma \lemB }
Let $\alpha\in\Cal D$ represent $\varphi(\tau)$.
Let $C$ be a component of the CRS for $\alpha$. Then
$C$ cannot separate $i$ and $j$ for $3\le i<j\le n$.
\endproclaim

\demo{ Proof }
Let $\beta\in\Cal D$ represent $\varphi(\sigma_{ij}^2\sigma_1^{-2})$.
Since $\alpha$ and $\beta$ commute, $\beta$ preserves $C$ by Proposition \propCRS(b).
Hence $C$ cannot separate $i$ and $j$
by Lemma \lemA\ applied to $\beta$ (note that $\beta$ is conjugate to $\sigma_1^2\sigma_3^{-2}$; see
the beginning of \S\sectProof.2).
\qed\enddemo

\proclaim{ Lemma \lemC }
Let $\alpha\in\Cal D$ represent $\varphi(\sigma_1^2\sigma_3^{-2})$.
Suppose that $n\ge 6$.
Let $C$ be a component of the CRS for $\alpha$. Then:

\smallskip\noindent
(a). $C$ cannot separate $1$ and $2$. It cannot separate $3$ and $4$.

\smallskip\noindent
(b). $C$ cannot separate $i$ and $j$ for $5\le i<j\le n$.

\smallskip\noindent
(c). $C$ cannot separate $\{1,2,3,4\}$ from $\{5,\dots,n\}$.

\smallskip\noindent
(d). $C$ cannot encircle $5,\dots,n$.
\endproclaim

\demo{ Proof }
(a). Follows from Lemma \lemA\ and Lemma \lemInvar.

\smallskip
(b). Let $\beta\in\Cal D$ represent $\varphi(\sigma_{ij}^{2}\sigma_1^{-2})$.
Since $\alpha$ and $\beta$ commute, $\beta$ preserves $C$ by Proposition \propCRS(b).
Hence $C$ cannot separate $i$ and $j$
by Lemma \lemA\ applied to $\beta$ (see the proof of Lemma \lemB).

\smallskip
(c). Suppose that $C$ separates $1,2,3,4$ from $5,6,\dots,n$.
Let $\beta\in\Cal D$ represent $\varphi(\sigma_1^{2}\sigma_5^{-2})$.
Then $\beta$ is conjugate to $\alpha$. Let $\gamma\in\Cal D$ be a conjugating element.
Then $\gamma(C)$ is a component of the CRS for $\beta$ and it separates the punctures
$1,2,5,6$ from all the other punctures.
Since $\alpha$ and $\beta$ commute, $\beta$ preserves $C$. This is impossible because
the geometric intersection number of $C$ and $\gamma(C)$ is nonzero.

\smallskip
(d). Combine (a), (c), and Lemma \lemInvar.
\qed\enddemo

\proclaim{ Lemma \lemD }
Let $\alpha\in\Cal D$ represent $\varphi(\sigma_1^2\sigma_3^{-2})$.
Suppose that $\alpha$ is reducible non-periodic.
Then the CRS for $\alpha$ has exactly two components:
one of them encircles $1$ and $2$, and the other one encircles $3$ and $4$.
\endproclaim

\demo{ Proof } 
If $n\ge 6$, the result follows from Lemma \lemInvar\ and Lemma \lemC.
Suppose that $n=5$ and the CRS is not as stated.
By combining Lemma \lemInvar\ with
Lemma \lemDegPer, we conclude that the CRS consists of a single circle which encircles 1,2,3,4.
The interior braid cannot be periodic by (\eqLK), hence it is pseudo-Anosov.
Therefore, $Z(\alpha;\P_5)\cong\Z^2$ by Proposition \propZ(b) whence
$Z(\alpha;\J_5)=\Z$. This contradicts Lemma \lemZgen.
\qed\enddemo

\proclaim{ Lemma \lemE }
$\varphi(\sigma_1\sigma_3^{-1})$ is conjugate in $\B_n$ to $\sigma_1\sigma_3^{-1}$.
\endproclaim

\demo{ Proof }
Let $\alpha\in\Cal D$ represent $\varphi(\sigma_1^2\sigma_3^{-2})$.
If $\alpha$ is pseudo-Anosov, then $Z(\alpha;\P_n)\cong\Z^2$ by Proposition \propZ(b), hence
$Z(\alpha;\J_n)$ is abelian which contradicts Lemma \lemZgen.
If $\alpha$ is periodic, then it is
a power of $\Delta^2$ by Proposition \propPerPur.
This contradicts (\eqLK), hence $\alpha$ is reducible non-periodic and
its CRS is as stated in Lemma \lemD.


Suppose that $\hat\alpha$ is pseudo-Anosov.
Then $Z(\hat\alpha;\P_n)\cong\Z^2$ by Proposition \propZ(b) whence
$Z(\alpha;\P_n)\cong\Z^4$ by Proposition \propZ(c) and therefore $Z(\alpha;\J_n)$ is abelian which
contradicts Lemma \lemZgen.
Thus, $\hat\alpha$ is periodic. By Proposition \propPerPur\ this means that $\hat\alpha$ is
a power of $\Delta^2$. This fact combined with (\eqLK) implies $\hat\alpha=1$.
It follows that $\varphi(\sigma_1^2\sigma_3^{-2})$ is conjugate to
$\sigma_1^{2k}\sigma_3^{-2k}$ for some $k$, and we have $k=1$ by (\eqLK).
The uniqueness of roots up to conjugation [\refGM] implies
that $\varphi(\sigma_1\sigma_3^{-1})$ is conjugate to $\sigma_1\sigma_3^{-1}$.
\qed\enddemo


\proclaim{ Lemma \lemF }
$\varphi(\tau)$ is conjugate in $\P_n$ to $\tau$.
\endproclaim

\demo{ Proof }
%
%
Let $\alpha\in\Cal D$ represent $\varphi(\tau)$.
By Proposition \propCRS, it cannot be pseudo-Anosov because
it commutes with $\varphi(\sigma_1\sigma_3^{-1})$ which is reducible non-periodic by Lemma \lemE.
If $\alpha$ were periodic, then it would be
a power of $\Delta^2$ by Proposition \propPerPur.
This contradicts (\eqLK), hence $\alpha$ is reducible.

Let $C$ be the CRS for $\alpha$. 
By Lemmas \lemA\ and \lemB, one of the following three cases takes place.

\smallskip
Case 1. $C$ is connected, the punctures $1$ and $2$ are inside $C$, all the other punctures are
outside $C$.
Then the tubular braid $\hat\alpha$ cannot be pseudo-Anosov because $\alpha$ commutes with
$\varphi(\sigma_1\sigma_3^{-1})$, hence it preserves a circle which separates $3$ and $4$ from
$5,\dots,n$. Hence $\hat\alpha$ is periodic which contradicts (\eqLK) combined with
Proposition \propPerPur. Thus this case is impossible.

\smallskip
Case 2. $C$ is connected, the punctures $1$ and $2$ are outside $C$, all the other punctures are
inside $C$. This case is also impossible and the proof is almost the same as in Case 1.
To show that $\hat\alpha$ cannot be pseudo-Anosov,
we note that $\alpha$ preserves a curve which encircles only $1$ and $2$.

\smallskip
Case 3. $C$ has two components: $c_1$ and $c_2$ which encircle $\{1,2\}$ and
$\{3,\dots,n\}$ respectively.
The interior braid $\alpha_{[2]}$ cannot be pseudo-Anosov by the same reasons as in Case 1:
because $\alpha$ preserves a circle separating $3$ and $4$ from $5,\dots,n$.
Hence $\alpha_{[2]}$ is periodic. Using (\eqLK), we conclude that $\alpha$ is a conjugate of
$\tau$.
Since the elements of $Z(\tau;\B_n)$ realize any permutation of $\{1,2\}$ and $\{3,\dots,n\}$,
the conjugating element can be chosen in $\P_n$.
\qed\enddemo

\proclaim{ Lemma \lemG }
There exists $\gamma\in\P_n$ such that 
$$
   \varphi(\sigma_1\sigma_i^{-1})=(\sigma_1^{-1}\sigma_i)^\gamma
   \quad\text{for $i=3,\dots,n$}.                                           \eqno(\eqLemG)
$$
\endproclaim

\demo{ Proof }
Due to Lemma \lemF, without loss of generality we may assume that $\varphi(\tau)=\tau$ and
$\tau(C)=C$ where $C$ is the CRS for $\tau$ consisting of two round circles $c_1$ and $c_2$
which encircle $\{1,2\}$ and $\{3,\dots,n\}$ respectively.

By Lemma \lemNonPure,
$\psi_{2,n-2}$ restricts to an isomorphism
$\psi:\P_2\times\B_2\times\B_{n-2}\to Z(\tau):=Z(\tau;\B_n)$.
Let $\pi_1:Z(\tau)\to\P_2$
and $\pi_3:Z(\tau)\to\B_{n-2}$ be defined as
$\pi_i=\operatorname{pr}_i\circ\psi^{-1}$.

Let $H=\pi_1^{-1}(1)\cap\B'_n$; note that the elements of $\pi_1^{-1}(1)$ correspond to geometric braids
whose first two strings are inside the cylinder $[0,1]\times c_1$ and the other strings are inside the cylinder
$[0,1]\times c_2$.
Then $\pi_3|_H:H\to\B_{n-2}$ is an isomorphism and its inverse is given by
$Y\mapsto\psi_{2,n-2}(1,\sigma_1^{-e(X)},Y)$, that is 
$\sigma_i\mapsto\sigma_1^{-1}\sigma_{i+2}$,
$i=1,\dots,n-3$.

Let us show that $\varphi(H)=H$. Indeed, let $X\in H$. Since $X\in Z(\tau;\B'_n)$ and $\varphi(\tau)=\tau$,
we have $\varphi(X)\in Z(\tau;\B'_n)$. The fact that $\pi_1(X)=1$ follows from (\eqLK)
applied to a power of $X$ belonging to $\J_n$. Hence $\varphi(H)\subset H$. By the same arguments
$\varphi^{-1}(H)\subset H$.

Thus $\varphi|_H$ is an automorphism of $H$ and we have $H\cong\B_{n-2}$. Hence,
by Dyer and Grossman's result [\refDG] cited after the statement of Theorem 1,
there exists $\gamma\in H$ such that $\ad\gamma\varphi|_H$ is either $\id_H$ or $\Lambda|_H$.
The latter case is impossible by (\eqLK). Thus there exists $\gamma\in\B_n$ such that
(\eqLemG) holds.

It remains to show that $\gamma$ can be chosen in $\P_n$. By replacing $\gamma$
with $\sigma_1\gamma$ if necessary, we may assume that $1$ and $2$ are fixed by $\gamma$.
By combining (\eqLKgamma), (\eqLK), and (\eqLemG), we conclude that $\gamma(\{i,j\})=\{i,j\}$
for any $i,j\in\{3,\dots,n\}$ and the result follows.
\qed\enddemo


\subhead\sectProof.2. Conjugates of $\sigma_1$ and simple curves which connect punctures
\endsubhead
We fix $n\ge2$ and we consider $\Bbb D$ and the set of punctures $X_n=\{1,\dots,n\}\subset\Bbb D$ as above.
Let $\Cal I$ be the set of all smooth simple curves (embedded segments) $I\subset\Bbb D$ such that
$\partial I\subset X_n$ and $I^\circ\subset\Bbb D\setminus X_n$. Here we denote
$I^\circ=I\setminus\partial I$ and $\partial I=\{a,b\}$ where $a$ and $b$ are the ends of $I$.
Recall that we write $I\sim I_1$ if $I_1=\alpha(I)$ for some $\alpha\in\Cal D_0$ (see \S\sectNT),
i.~e., if $I$ and $I_1$ belong to the same connected component of $\Cal I$.

Let $I\in\Cal I$ and
let $\beta\in\Cal D$ be such that $\beta(I)$ is the straight line segment $[1,2]$.
Then we define the braid $\sigma_{I}$ as $\sigma_1^\beta$.
It is easy to see that $\sigma_I$ depends only on the connected component of $\Cal I$ that contains $I$.
The CRS for $\sigma_I$ is a single closed curve which encloses $I$ and separates it from $X_n\setminus\partial I$.
By definition, all conjugates of $\sigma_1$ are obtained in this way.
In particular, we have $\sigma_i=\sigma_{[i,i+1]}$ and $\sigma_{ij}=\sigma_I$ for an
embedded segment $I$  which connects $i$ to $j$ passing through
the upper half-plane.

\proclaim{ Lemma \lemConjI }
For any $\beta\in\Cal D$, $I\in\Cal I$, we have
$\sigma_{\beta(I)}^\beta = \sigma_I$.
\qed
\endproclaim

In this notation, a corollary of Lemma \lemG\ can be formulated as follows.

\proclaim{ Lemma \lemH } Let $n\ge5$ and let $\varphi\in\Aut(B'_n)$ be as in \S\sectProof.1.

\smallskip\noindent
(a). Let $I,J\in\Cal I$ be such that
$\Card(I\cap J)=\Card(\partial I\cap\partial J)=1$
{\rm(}i.e., $I\cap J$ is a common endpoint of $I$ and $J${\rm)}.
Then there exist 
$I_1,J_1\in\Cal I$ such that $I_1\cup J_1$ is homeomorphic to $I\cup J$ and
$$
  \varphi(\sigma_I\sigma_J^{-1})=\sigma_{I_1}\sigma_{J_1}^{-1}, \qquad
  \varphi(\sigma_I^{-1}\sigma_J)=\sigma_{I_1}^{-1}\sigma_{J_1}.                        \eqno(\eqLemHa)
$$
(b). Let $I,J\in\Cal I$, $I\cap J=\varnothing$. Then the conclusion is the same as in Part (a).

\smallskip\noindent
(c). Let $I$ and $J$ be as in Part (a) and let $K\in\Cal I$ be such that $K\cap(I\cup J)=\varnothing$.
Then there exist 
$I_1,J_1,K_1\in\Cal I$ such that
$I_1\cup J_1\cup K_1$ is homeomorphic to $I\cup J\cup K$, and (\eqLemHa) holds as well as
$$
  \varphi(\sigma_K\sigma_I^{-1})=\sigma_{K_1}\sigma_{I_1}^{-1},\qquad
  \varphi(\sigma_K\sigma_J^{-1})=\sigma_{K_1}\sigma_{J_1}^{-1}.                        \eqno(\eqLemHc)
$$
\endproclaim

\demo{ Proof } (c).
Let $\gamma$ be as in Lemma \lemG\ and
let $\beta\in\Cal D$ be such that $\beta(K)=[1,2]$, $\beta(I)=[3,4]$, and $\beta(J)=[4,5]$. 
We set $K_1=\alpha^{-1}(K)$, $I_1=\alpha^{-1}(I)$, $J_1=\alpha^{-1}(J)$ where
$\alpha=\beta^{-1}\gamma\varphi(\beta)$.
Then we have
$$
\xalignat2
   \varphi(\sigma_K\sigma_I^{-1})
     &= \varphi((\sigma_1\sigma_3^{-1})^\beta) &&\text{by definition of $\sigma_I$ and $\sigma_K$}\\
     &= (\sigma_1\sigma_3^{-1})^{\gamma\varphi(\beta)}    &&\text{by Lemma \lemG}\\
     &= \sigma_{K_1}\sigma_{I_1}^{-1}            &&\text{by Lemma \lemConjI}
\endxalignat
$$
and, similarly, $\varphi(\sigma_K\sigma_J^{-1})=\sigma_{K_1}\sigma_{J_1}^{-1}$.
Since $\sigma_K$ commutes with $\sigma_I$ and $\sigma_J$,
we have
$\sigma_I^\eps\sigma_J^{-\eps} = (\sigma_K\sigma_I^{-1})^{-\eps}(\sigma_K\sigma_J^{-1})^\eps$, $\eps=\pm1$,
thus (\eqLemHc) implies (\eqLemHa).


\smallskip
(a). 
Since $\card(\partial I\cup\partial J)=3$ and $n\ge5$, we can choose $K\in\Cal I$
disjoint from $I\cup J$ (which is an embedded segment, hence its complement is
connected) and the result follows from (c).

\smallskip
(b). The same proof as for Part (c) but with $\beta(I)=[1,2]$ and $\beta(J)=[3,4]$.
\qed\enddemo

\proclaim{ Lemma \lemUniq }
Let $I,J,I_1,J_1\in\Cal I$ be such that $I\cap J=I_1\cap J_1=\varnothing$.
Suppose that $\sigma_I\sigma_J^{-1} = \sigma_{I_1}\sigma_{J_1}^{-1}$.
Then $I\sim I_1$ and $J\sim J_1$.
\endproclaim

\demo{ Proof } It is enough to observe that the CRS for $\sigma_I\sigma_J^{-1}$ is the $\partial U_I\cup\partial U_J$
where $U_I$ and $U_J$ are $\varepsilon$-neighbourhoods of $I$ and $J$ for $0<\varepsilon\ll 1$
(this fact follows, for example, from Lemma \lemD\ and Proposition \propConjCRS).
\qed
\enddemo

Note that when $[\sigma_I,\sigma_J]\ne1$, the statement of Lemma \lemUniq\
is wrong. Indeed, in this case by Lemma \lemConjI\ we have
$\sigma_I\sigma_J^{-1}=\sigma_{\gamma(I)}\sigma_{\gamma(J)}^{-1}$ for $\gamma=\sigma_I\sigma_J^{-1}$
whereas $\sigma_I\ne\sigma_{\gamma(I)}$ and $\sigma_J\ne\sigma_{\gamma(J)}$.

\smallskip

Given $I,J\in\Cal I$, the {\it geometric intersection number} $I\cdot J$ of $I$ and $J$ is defined
as the minimum of the number of intersection points of $I_1^\circ$ and $J_1^\circ$
over all pairs $(I_1,J_1)\in\Cal I^2$ such that $I\sim I_1$, $J\sim J_1$, and
$I_1$ is transverse to $J_1$. In this case we say that $I_1$ and $J_1$ {\it realize} $I\cdot J$.

%

\midinsert
\centerline{\epsfxsize=90mm\epsfbox{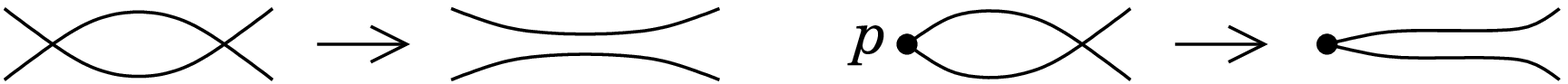}}
\botcaption{ Figure \figMove } Digon removal ($p$ is a puncture) \endcaption
\endinsert

If $I,J\in\Cal I$ are transverse to each other, we say that a closed embedded disk $D$ is a
{\it digon between $I$ and $J$} if $D$ is the closure of a component of $\Bbb D\setminus(I\cup J)$,
and $\partial D$ is a union of an arc of $I$ and an arc of $J$. The common ends of these arcs
are called the {\it corners} of $D$. We say that $(I',J')$ is obtained from $(I,J)$ by a
{\it digon removal} if it is obtained by one of the modifications in Figure \figMove\
performed in a neighbourhood of a digon between $I$ and $J$ one of whose corners is not in $X_n$.
The inverse operation is called a {\it digon insertion}.

The following two lemmas have a lot of analogs in the literature but it is easier
to write (and to read) a proof than to search for an appropriate reference.

\proclaim{ Lemma \lemIone } Let $I,J\in\Cal I$ be transverse to each other.
Then a pair of segments realizing $I\cdot J$ can be obtained from $(I,J)$ by successive
digon removals.
\endproclaim

\demo{ Proof }
Isotopies of $I$ and of $J$ which transform $(I,J)$ to a pair of segments realizing $I\cdot J$ can
be perturbed into a sequence of digon removals and digon insertions. So, it is enough to prove
the following ``diamond lemma": if $(I_1,J_1)$ and $(I_2,J_2)$ are obtained from $(I,J)$ by two
different digon removals, then either the pair $(I_1\cup J_1,I_1)$ is isotopic to $(I_2\cup J_2,I_2)$,
or $(I_1,J_1)$ and $(I_2,J_2)$ admit digon removals with the same result. We leave it to the reader to
check this statement (see Figure \figDiamond).
\qed\enddemo

\midinsert
\centerline{\epsfxsize=70mm\epsfbox{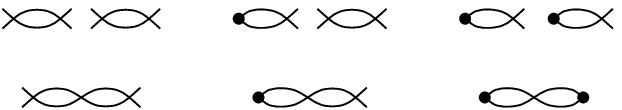}}
\botcaption{ Figure \figDiamond } Cases to consider in the diamond lemma \endcaption
\endinsert

\proclaim{ Lemma \lemItwo } Let $I_1,\dots,I_m\in\Cal I$. Then there exist
$I'_1,\dots,I'_m\in\Cal I$ such that $I_i\sim I'_i$ for any $i=1,\dots,m$, and
$(I'_i,I'_j)$ realizes $I_i\cdot I_j$ for any distinct $i,j=1,\dots,m$.
\endproclaim

\demo{ Proof } Induction on the total number of intersection points.
If $(I_i,I_j)$ does not realize $I_i\cdot I_j$, then by Lemma \lemIone\
there is a digon $D$ between $I_i$ and $I_j$.
We can remove $D$ so that the union of all segments is modified only near the corners of
$D$. Then the total number of intersection points strictly decreases.
\qed\enddemo


\subhead\sectProof.3. End of the proof
\endsubhead
Now we are ready to complete the proof of Theorem 1 for $n\ge 5$.
Let $n\ge5$ and let $\varphi$ be an automorphism of $B'_n$.
By [\refLin; Theorem C], we may assume that either $\mu'\varphi=\mu'$, or $n=6$ and $\mu'\varphi=\nu\mu'$
where $\nu$ is as in \S\sectRepr. However, $\mu'\varphi\ne\nu\mu'$ by Lemma \lemReprSix.
So, we assume that $\mu'\varphi=\mu'$. Then Lemma \lemRepr\ implies that the automorphism $\varphi_*$
of $\J_n^\ab$ induced by $\varphi$ is $\pm\id$. By composing $\varphi$ with $\Lambda$ if necessary,
we may assume that $\varphi_*=\id$
(recall that $\Lambda$ is the automorphism of $\B_n$ which takes each $\sigma_i$ to $\sigma_i^{-1}$).
By Lemma \lemG\ we may also assume that
$$
   \varphi(\sigma_1\sigma_i^{-1})=\sigma_1\sigma_i^{-1}
   \quad\text{for all $i=3,\dots,n-1$}                                                  \eqno(\eqA)
$$
(otherwise we compose $\varphi$ with $\ad\gamma$ for the $\gamma$ from Lemma \lemG). Hence
$$
    \varphi(\sigma_i\sigma_j^{-1}) = \sigma_i\sigma_j^{-1}\;\text{ and }\;
    \varphi(\sigma_i^{-1}\sigma_j) = \sigma_i^{-1}\sigma_j
    \quad\text{for all $i,j\in\{3,\dots,n-1\}$.}  \eqno(\eqB)
$$
Indeed, $\sigma_i\sigma_j^{-1}=(\sigma_1\sigma_i^{-1})^{-1}(\sigma_1\sigma_j^{-1})$ and
$\sigma_i^{-1}\sigma_j=(\sigma_1\sigma_i^{-1})(\sigma_1\sigma_j^{-1})^{-1}$.

Let $I_1,J_1,K_1\in\Cal I$ be as in Lemma \lemH(c) where we set $I=[1,2]$, $J=[2,3]$, and $K=[4,5]$.
By combining (\eqLemHc) 
with (\eqA) for $i=4$, we obtain $\sigma_1\sigma_4^{-1} = \sigma_{I_1}\sigma_{K_1}^{-1}$. Hence
$I_1\sim[1,2]$ and $K_1\sim[4,5]$ by Lemma \lemUniq.
Thus, if we set $L=J_1$, then (\eqLemHa) reads as
$$
     \varphi(\sigma_1\sigma_2^{-1}) = \sigma_{1}\sigma_{L}^{-1}\quad\text{ and }\quad
     \varphi(\sigma_1^{-1}\sigma_2) = \sigma_{1}^{-1}\sigma_{L}.                          \eqno(\eqC)
$$
By (\eqLK) for $\sigma_2\sigma_4^{-2}$ and by the claim $I_1\cup J_1\cong I\cup J$ of Lemma \lemH(c) we also have
$$
     \partial L=\{2,3\} \quad\text{ and }\quad  L\cdot[1,2]=0.                             \eqno(\eqD)
$$
By combining (\eqA) with (\eqC), we obtain
$$
     \varphi(\sigma_{i}\sigma_2^{-1}) = \sigma_{i}\sigma_{L}^{-1}\;\;\text{ and }\;\;
     \varphi(\sigma_{i}^{-1}\sigma_2) = \sigma_{i}^{-1}\sigma_{L}
                  \quad\text{for all $i=3,\dots,n-1$.}                                     \eqno(\eqE)
$$
This fact combined with Lemma \lemH(b) and Lemma \lemUniq\ implies
$$
     L\cdot[i,i+1]=0\qquad\text{ for all $i=4,\dots,n-1$}.                                 \eqno(\eqF)
$$
Indeed, for any $i=4,\dots,n-1$, by Lemma \lemH(b) we have
$\varphi(\sigma_2\sigma_i^{-1})=\sigma_{I_1}\sigma_{J_1}^{-1}$ for some disjoint $I_1,J_1\in\Cal I$.
On the other hand, $\varphi(\sigma_2\sigma_i^{-1})=\sigma_L\sigma_i^{-1}$ by (\eqE).
Hence $I_1\sim L$ and $J_1\sim[i,i+1]$ by Lemma \lemUniq\ whence (\eqF) because $I_1\cap J_1=\varnothing$.

\proclaim{ Lemma \lemL } $L\cdot[3,4]=0$.
\endproclaim

\demo{ Proof }
By Lemma \lemH(a) applied to $I=[2,3]$ and $J=[3,4]$,
there exist $I_1,J_1$ such that $I_1\cup J_1\cong[2,4]$ (thus $I_1\cdot J_1=0$)
and $\varphi(\sigma_2\sigma_3^{-1})=\sigma_{I_1}\sigma_{J_1}^{-1}$.
By combining this fact with (\eqE) for $i=3$,
we obtain $\sigma_L\sigma_3^{-1}=\sigma_{I_1}\sigma_{J_1}^{-1}$. Hence, by Lemma \lemGarsideXY,
there exists $\gamma\in\Cal D$ such that
$\sigma_L^\gamma=\sigma_{I_1}$ and $\sigma_3^\gamma = \sigma_{J_1}$
whence $\gamma(I_1)\sim L$ and $\gamma(J_1)\sim[2,3]$ by Lemma \lemConjI. Thus
$L\cdot[3,4]=I_1\cdot J_1=0$.
\qed\enddemo

\demo{ Proof of Lemma \lemL\ for $n\ge 6$ not using Garside theory}
Let $n\ge 6$. We apply the same arguments that we used to obtain (\eqC)--(\eqE) but
we set here $I=[3,4]$, $J=[2,3]$, $K=[5,6]$. So, let
$I_1,J_1,K_1$ be as in Lemma \lemH(c) for the given choice of $I,J,K$.
By combining (\eqLemHc) with (\eqB) and (\eqE), we obtain
$\sigma_{I_1}\sigma_{K_1}^{-1}=\sigma_3\sigma_5^{-1}$ and
$\sigma_{J_1}\sigma_{K_1}^{-1}=\sigma_L\sigma_5^{-1}$. Then Lemma \lemUniq\
yields $I_1=[3,4]$, $J_1=L$, and $K_1=[5,6]$.
By Lemma \lemH(c), $I_1\cup J_1$ is homeomorphic to $I\cup J$, hence
$L\cdot[3,4]=I\cdot J=0$.
\qed\enddemo

Further, Lemma \lemL\ combined with (\eqD) and (\eqF), yields $L\cdot[i,i+1]=0$ for any
$i\in\{1\}\cup\{3,\dots,n-1\}$. 
By Lemma \lemItwo\ this implies that
$L\sim L_1$ where $L_1\in\Cal I$ is such that $[1,2]\cup L_1\cup[3,n]$ is homeomorphic to a segment.
Hence, up to composing $\varphi$ with $\tilde\beta$ where
$\beta\in\Cal D$, $\beta([1,n])=[1,2]\cup L_1\cup[3,n]$,
we may assume that $\sigma_L=\sigma_2$ in (\eqA)--(\eqC) and (\eqE). This means
that $\varphi(\sigma_i^\eps\sigma_j^{-\eps})=\sigma_i^\eps\sigma_j^{-\eps}$
for any $\eps=\pm1$ and any $i,j\in\{1,\dots,n-1\}$.
To complete the proof of Theorem 1 for $n\ge 5$, it remains to note that the elements
$\sigma_i^{\eps}\sigma_j^{-\eps}$,  
$\eps=\pm1$,
generate $\B'_n$. Indeed, it is shown in [\refGL]
(see also [\refLin; \S1.8]) that $\B'_n$ is generated by
$u=\sigma_2\sigma_1^{-1}$, $v=\sigma_1\sigma_2\sigma_1^{-2}=(\sigma_2^{-1}\sigma_1)(\sigma_2\sigma_1^{-1})$,
$w=(\sigma_2\sigma_1^{-1})(\sigma_3\sigma_2^{-1})$, and
$c_i=\sigma_i\sigma_1^{-1}$, $i=3,\dots,n-1$.

\medskip\noindent
{\bf Remark \remProof.}
Our proof of Theorem 1 for $n\ge5$ essentially uses Lemma \lemG\ which is based on
Dyer-Grossman's result [\refDG] about $\Aut(B_n)$.
If $n\ge 6$, Lemma \lemG\ 
can be replaced by Lemma \lemGarside\ (see below).


\head\sectBfour. The case $n=4$
\endhead

Recall that $\B'_3$ is freely generated by $u=\sigma_2\sigma_1^{-1}$ and $t=\sigma_1^{-1}\sigma_2$
(see the Introduction). The group $\B'_4$ was computed in [\refGL], namely
$\B'_4=\bold K_4\rtimes\B'_3$ where
$\bold K_4$ is the kernel of the homomorphism $\B_4\to\B_3$, $\sigma_1,\sigma_3\mapsto\sigma_1$, $\sigma_2\mapsto\sigma_2$.
The group $\bold K_4$ is freely generated by $c=\sigma_3\sigma_1^{-1}$ and $w=\sigma_2c\,\sigma_2^{-1}$.
The action of $\B'_3$ on $\bold K_4$ by conjugation is given by
$$
  ucu^{-1}=w,\qquad uwu^{-1}=w^2c^{-1}w,\qquad tct^{-1}=cw,\qquad twt^{-1}=cw^2.           \eqno(\eqGL)
$$
Besides the elements $c,w,u,t$ of $\B'_4$, we consider also
$$
   d=\psi_{2,2}(\sigma_1^{-1};\sigma_1^2,\sigma_1^2)=\sigma_1^3\sigma_3^3\Delta^{-1}.
$$

\proclaim{ Lemma \lemZd }

\smallskip\noindent
(a). $Z(d^2;\B'_4)$ is a semidirect product of infinite cyclic groups
$\langle c\rangle\rtimes\langle d\rangle$ where
$d$ acts on $\langle c\rangle$ by $dcd^{-1}=c^{-1}$.

\smallskip\noindent
(b). $\langle c\rangle$ is a characteristic subgroup of $Z(d^2;\B'_4)$.
\endproclaim

\demo{ Proof } Let $G=Z(d^2;\B'_4)$.

(a). We have $G=Z(d^2;\B_4)\cap\ker e$ and, by [\refGMW; \S5],  $Z(d^2;\B_4)$ is the semidirect product
$\langle\sigma_1,\sigma_3\rangle\rtimes\langle d\rangle$ where
$d$ acts on $\langle\sigma_1,\sigma_3\rangle$ by $\sigma_1^d=\sigma_3$, $\sigma_3^d=\sigma_1$.

(b).
Let $x$ be the image of $c$ by an automorphism of $G$. Then (a) implies that
$x$ generates a normal subgroup of $G$ and $x$ is not a power of another element of $G$.
It follows that $x\in\{c,c^{-1}\}$.
\qed\enddemo

\proclaim{ Lemma \lemZall } All conjugacy classes of $\B_4$ which are contained in $\B'_4$ are presented in Table 1.
The corresponding centralizers are isomorphic to the groups indicated in this table.
\endproclaim

\demo{ Proof }
First, note that $\B'_n$ is normal in $\B_n$, hence for $X\in\B'_n$, the centralizer $Z(X;\B'_n)$
depends only on the conjugacy class of $X$ in $\B_n$
(though this class may split into several classes in $\B'_n$).

The centralizers in $\B_4$ can be computed by a straightforward application
of 
[\refGMW; Propositions 4.1] and Proposition \propCRS.
In the computation of $Z(\Delta^{2k+1}\sigma_2^{-12k-6})$ (which is, by the way, generated by
$\Delta$ and $\sigma_2$), we use the fact that $Z(\Delta_3^{2k+1};\B_3)=\langle\Delta\rangle$.
This fact can be derived either from [\refGMW; Proposition 3.5] or from the uniqueness of Garside
normal form in $\B_3$.

In all the cases except, maybe, the following two ones, the computation of
$Z(X;\B'_4)$ is evident.

\smallskip
1). $X=Y^k$, $Y=\psi_{3,1}(\sigma_1^{-2};\Delta_3^{2},1)$, $k\ne0$. 
The group $Z(X;\B_4)$ is generated by $\psi_{3,1}(\sigma_1^2;1,1)$, $\sigma_1$, and $\sigma_2$.
Since $\psi_{3,1}(\sigma_1^2;1,1)=Y\Delta_3^{2}$ and
$\Delta_3\in\B_3$, we can choose
$Y,\sigma_1,\sigma_2$ for a generating set, and the result follows because $e(Y)=0$.

\smallskip
2).
$X=\Delta^{2k}\sigma_1^{-6k}$, $k\ne 0$. 
The group $Z(X;\B_4)$ is the isomorphic image of $\B_{2,1}\times\Z$
under the mapping $f:(X,m)\mapsto\psi_{2,1}(X;\sigma_1^m,1)$.
Hence $Z(X;\B'_4)$ is the isomorphic image of $\B_{2,1}$
under the mapping $X\mapsto f(X,-e(f(X,0))$.
\qed\enddemo

\midinsert
\noindent Table 1. Centralizers of elements of $\B'_4$;

\noindent white/grey region $\Rightarrow$ the associated braid is periodic/pseudo-Anosov;

\noindent in $Z(d^{2k};\B_4)$ we mean $f:\Z\to\Aut(\Z\times\Z)$, $f(1)(x,y)=(y,x)$.
\smallskip
\vbox{\offinterlineskip
\def\e#1{\hskip-2pt\epsfxsize=9mm\lower3pt\hbox{\epsfbox{#1.eps}}}
\def\TA{$\psi_{3,1}(\sigma_1^{-2k};\Delta_3^{2k},1)$, $k\ne0$}

\def\TC{$d^{2k}$, $k\ne0$}
\def\TD{$d^{2k}c^l$, $l\not\in\{0,\pm 6k\}$}
\def\TE{$d^{2k+1}$}
\def\TF{$\Delta^{2k}\sigma_1^{-12k}$, $k\ne0$}
\def\TG{$\Delta^{2k+1}\sigma_2^{-12k-6}$}

\hrule
\halign{&\vrule#&\strut\quad\;#\hfill\cr
height3pt&\omit&\omit& \omit  &\omit& \omit              &\omit& \omit           &\omit& \omit              &\cr
& CRS          &\omit& $X$    &\omit& $Z(X;\B_4)$        &\omit& $Z(X;\B'_4)$    &\omit& $Z(X;\B'_4)^\ab\;$ &\cr
height3pt&\omit&\omit& \omit  &\omit& \omit              &\omit&  \omit          &\omit& \omit              &\cr
\noalign{\hrule}
height3pt&\omit&\omit& \omit  &\omit& \omit              &\omit&  \omit          &\omit& \omit              &\cr
&\e{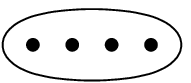}        &\omit& $1$    &\omit& $\B_4$             &\omit& $\B'_4$         &\omit& \omit              &\cr
height3pt&\omit&\omit& \omit  &\omit& \omit              &\omit&  \omit          &\omit& \omit              &\cr
\noalign{\hrule}
height3pt&\omit&\omit& \omit  &\omit& \omit              &\omit&  \omit          &\omit& \omit             &\cr
&\e{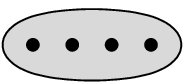}        &\omit& \omit  &\omit& $\Z^2$             &\omit& $\Z$            &\omit& \omit             &\cr
height3pt&\omit&\omit& \omit  &\omit& \omit              &\omit&  \omit          &\omit& \omit             &\cr
\noalign{\hrule}
height3pt&\omit&\omit& \omit  &\omit& \omit              &\omit&  \omit          &\omit& \omit             &\cr
&\e{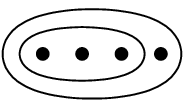}        &\omit& \TA    &\omit& $\B_3\times\Z$     &\omit& $\B'_3\times\Z$ &\omit& $\Z^3$            &\cr
height3pt&\omit&\omit& \omit  &\omit& \omit              &\omit&  \omit          &\omit& \omit             &\cr
\noalign{\hrule}
height3pt&\omit&\omit& \omit  &\omit& \omit              &\omit&  \omit          &\omit& \omit             &\cr
&\e{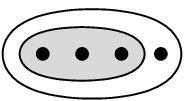}        &\omit& \omit  &\omit& $\Z^3$             &\omit& $\Z^2$          &\omit& \omit             &\cr
height3pt&\omit&\omit& \omit  &\omit& \omit              &\omit&  \omit          &\omit& \omit             &\cr
\noalign{\hrule}
height3pt&\omit&\omit& \omit  &\omit& \omit              &\omit&  \omit          &\omit& \omit             &\cr
&\e{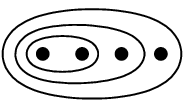}        &\omit& \omit  &\omit& $\Z^3$             &\omit& $\Z^2$          &\omit& \omit             &\cr
height3pt&\omit&\omit& \omit  &\omit& \omit              &\omit&  \omit          &\omit& \omit             &\cr
\noalign{\hrule}
height3pt&\omit&\omit& \omit  &\omit& \omit              &\omit&  \omit          &\omit& \omit             &\cr
height3pt&\omit&\omit& \omit  &\omit& \omit              &\omit&  \omit          &\omit& \omit             &\cr
&\omit         &\omit& \TC    &\omit& $\Z^2\rtimes_f\Z$  &\omit& $\Z\rtimes\Z$   &\omit& $\Z\times\Z_2$    &\cr
height3pt&\omit&\omit& \omit  &\omit& \omit              &\omit&  \omit          &\omit& \omit             &\cr
&\e{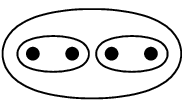}        &\omit& \TD    &\omit& $\Z^3$             &\omit& $\Z^2$          &\omit& \omit             &\cr
height3pt&\omit&\omit& \omit  &\omit& \omit              &\omit&  \omit          &\omit& \omit             &\cr
&\omit         &\omit& \TE    &\omit& $\Z^2$             &\omit& $\Z$            &\omit& \omit             &\cr
height3pt&\omit&\omit& \omit  &\omit& \omit              &\omit&  \omit          &\omit& \omit             &\cr
\noalign{\hrule}
height3pt&\omit&\omit& \omit  &\omit& \omit              &\omit&  \omit          &\omit& \omit             &\cr
&\e{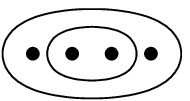}        &\omit& \TF    &\omit& $\B_{2,1}\times\Z$ &\omit& $\B_{2,1}$      &\omit& $\Z^2$            &\cr
height3pt&\omit&\omit& \omit  &\omit& \omit              &\omit&  \omit          &\omit& \omit             &\cr
&\omit         &\omit& \TG    &\omit& $\Z^2$             &\omit& $\Z$            &\omit& \omit             &\cr
height3pt&\omit&\omit& \omit  &\omit& \omit              &\omit&  \omit          &\omit& \omit             &\cr
\noalign{\hrule}
height3pt&\omit&\omit& \omit  &\omit& \omit              &\omit&  \omit          &\omit& \omit             &\cr
&\e{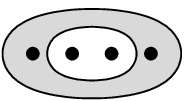}        &\omit& \omit  &\omit& $\Z^3$             &\omit& $\Z^2$          &\omit& \omit             &\cr
height3pt&\omit&\omit& \omit  &\omit& \omit              &\omit&  \omit          &\omit& \omit             &\cr
}\hrule}
\endinsert

Let $\varphi\in\Aut(\B'_4)$.

\proclaim{ Lemma \lemInvd } $\varphi(d)$ is conjugate in $\B_4$ to $d^{\pm1}$.
\endproclaim

\demo{ Proof } Let $x=\varphi(d^2)$. Since $Z(x;\B'_4)\cong Z(d^2;\B'_4)$, we see in Table 1 that
$\varphi(\langle d^2\rangle)\subset\langle d^2\rangle$.
By the same reasons we have
$\varphi^{-1}(\langle d^2\rangle)\subset\langle d^2\rangle$, thus
$\varphi(d^2)=d^{\pm2}$
and the result follows from the uniqueness of roots up to conjugation [\refGM].
\qed\enddemo

\proclaim{ Lemma \lemInvc } If $\varphi(d)=d$, then $\varphi(c)=c^{\pm1}$.
\endproclaim

\demo{ Proof } If $\varphi(d)=d$, then $\varphi(Z(d^2;\B'_4))=Z(d^2;\B'_4)$,
and we apply Lemma \lemZd.
\qed\enddemo

\proclaim{ Lemma \lemInvK } $\bold K_4$ is a characteristic subgroup in $\B'_4$.
\endproclaim

\demo{ Proof } Lemma \lemInvd\ combined with Lemma \lemInvc\ imply that $\varphi(c)$
is conjugate to $c$ in $\B_4$. Since $\bold K_4$ is the normal closure of $c$ in $\B_4$,
it follows that $\varphi(c)\in\bold K_4$. The same arguments can be applied to any other automorphism of $\B'_4$,
in particular, to $\varphi\ad\sigma_2$ whence $\varphi\ad\sigma_2(c)\in\bold K_4$. It remains to recall
that $\varphi\ad\sigma_2(c)=\varphi(w)$ and $\bold K_4=\langle c,w\rangle$.
\qed\enddemo

Let
$$
    S_1          =\left(\matrix 1 &-1\\ 0 & 1\endmatrix\right)\!,\,
    S_2          =\left(\matrix 1 & 0\\ 1 & 1\endmatrix\right)\!,\,
    T=S_1^{-1}S_2=\left(\matrix 2 & 1\\ 1 & 1\endmatrix\right)\!,\,
    U=S_2S_1^{-1}=\left(\matrix 1 & 1\\ 1 & 2\endmatrix\right)\!.
$$

\proclaim{ Lemma \lemTU } $T$ and $U$ generate a free subgroup of $\SL(2;\Z)$.
\endproclaim

\demo{ Proof } It is well known that the correspondence $\sigma_1\mapsto S_1$, $\sigma_2\mapsto S_2$
defines an isomorphism $\B_3/\langle\Delta^2\rangle\to\PSL(2,\Z)$, see, e.g., [\refMKS; \S3.5]
(this mapping is also a specialization of
the reduced Burau representation).
Since $u\mapsto U$ and $t\mapsto T$, the image of $\B'_3=\langle u,t\rangle$ is $\langle U,T\rangle$.
Hence $\langle U,T\rangle$ is free.
\qed\enddemo

\proclaim{ Lemma \lemK } If $\varphi|_{\bold K_4}=\id$, then $\varphi=\id$.
\endproclaim

\demo{ Proof } Let $\varphi|_{\bold K_4}=\id$. Since $\B'_4=\bold K_4\rtimes\B'_3$, we may write
$\varphi(u)=u_1 a$ and $\varphi(t)=t_1 b$ with $u_1,t_1\in\B'_3$ and $a,b\in\bold K_4$.
For $x\in\bold K_4$, we have $\varphi\tilde u(x)=\varphi(uxu^{-1})=u_1axa^{-1}u_1^{-1}=\tilde u_1\tilde a(x)$.
Since $\tilde u(x)\in\bold K_4$ and $\varphi|_{\bold K_4}=\id$, we conclude that
$\tilde u(x)=\varphi\ad u(x)=\tilde u_1\tilde a(x)$.
Similarly, $\tilde t(x)=\tilde t_1\tilde b(x)$. Thus,
$$
     \tilde u\,|_{\bold K_4}=\tilde u_1\tilde b\,|_{\bold K_4}
     \qquad\text{and}\qquad
     \tilde t\,|_{\bold K_4}=\tilde t_1\tilde a\,|_{\bold K_4}             \eqno(\eqUTK)
$$

Consider the homomorphism $\pi:\B'_4\to\Aut(\bold K_4^\ab)=\GL(2,\Z)$,
$x\mapsto(\tilde x)_*$;
here we identify  $\Aut(\bold K_4^\ab)$ with $\GL(2,\Z)$ by choosing the images of $c$ and $w$
as a base of $\bold K_4^\ab$.
It is clear that $\pi(a)=\pi(b)=1$ and it follows from (\eqGL) that
$\pi(tu^{-1})=T$ and $\pi(t)=U$. Thus, by Lemma \lemTU, the restriction of $\pi$ to $\B'_3$ is injective.
It follows from (\eqUTK) that $\pi(u_1)=\pi(u)$ and $\pi(t_1)=\pi(t)$.
Hence $u_1=u$ and $t_1=t$. Then it follows from (\eqUTK) that
$\tilde a\,|_{\bold K_4}=\tilde b\,|_{\bold K_4}=\id$.
Since $\bold K_4$ is free, its center is trivial, and we obtain $a=b=1$.
Thus $\varphi=\id$.
\qed\enddemo

\demo{ Proof of Theorem 1 for $n=4$ }
Let $\varphi\in\Aut(\B'_4)$.
By Lemma \lemInvd, we may assume that $\varphi(d)=d^{\pm1}$.
Then, by Lemma \lemInvc, we may assume that $\varphi(c)=c^{\pm1}$.
Since $c^\Delta=c^{-1}$, we may further assume that $\varphi(c)=c$.
By Lemma \lemInvK, $\varphi(c)$ and $\varphi(w)$ is a free base of $\bold K_4$.
Since $\varphi(c)=c$, it follows that $\varphi(w)=c^pw^{\pm1}c^q$, $p,q\in\Z$, see
[\refMKS; \S3.5, Problem 3].
We have
$\tilde\sigma_1(c)=c$ and
{\def\1{\bar 1}\def\2{\bar 2}\def\3{\bar 3}
$\tilde\sigma_1(w)= 1 2 3\1\2\1
                  = 1 2 3\2\1\2
                  = 1\3 2 3\1\2=c^{-1}w$}.
(here $1,\bar 1,2,\dots$ stand for $\sigma_1,\sigma_1^{-1},\sigma_2,\dots$).
Thus, by composing $\varphi$ with a power of $\tilde c$ and a power of $\tilde\sigma_1$ if necessary,
we may assume that $\varphi(w)=w^{\pm1}$.
For $\Phi=\Lambda\ad\sigma_1\ad\sigma_3\ad\Delta$,
we have $\Phi(c)=c$ and 
{\def\1{\bar 1}\def\2{\bar 2}\def\3{\bar 3}
$\Phi(w)=\1\3\2 3\1 2 1 3
        =\1 2\3\2 2 1\2 3 
        =\1 2 1 \3\2 3 
        = 2 1\2  2\3\2 = w^{-1}$}
hence,
by composing $\varphi$ with $\Phi$ if necessary,
we may assume that $\varphi(c)=c$ and $\varphi(w)=w$, thus
$\varphi|_{\bold K_4}=\id$ and the result
follows from Lemma \lemK. 
\qed\enddemo




\head\sectGarside. Appendix. Garside-theoretic lemmas
\endhead

Here, using Garside theory, we prove two statements one of which (Lemma \lemGarsideXY) is used only in the proof
of Theorem 1 for $n=5$, see the proofs of Lemma \lemL, and the other one (Lemma \lemGarside)
can be used in the proof of Theorem 1 for $n\ge6$ instead of Dyer-Grossman theorem, see Remark \remProof.

Let $n\ge3$ and
let $\Cal I$ and $\sigma_I\in\B_n$ for $I\in\Cal I$ be as in \S\sectI.

\proclaim{ Lemma \lemGarsideXY } Let $k,l\in\Z\setminus\{0\}$ and $I,J\in\Cal I$.
Suppose that $\sigma_I^k \sigma_J^l$ is conjugate to $\sigma_1^k\sigma_2^l$.
Then there exists $u\in\B_n$ such that
$\sigma_I^u = \sigma_1$ and $\sigma_J^u = \sigma_2$,
in particular, $I\cdot J=0$.
\endproclaim

\demo{ Proof } It follows from Corollary \corQP\ that there exists $u\in\B_n$ and
$p,q,r,s\in\Z\cap[1,n]$ such that $\sigma_I^u = \sigma_{pq}$ and $\sigma_J^u = \sigma_{rs}$,
This means that $\sigma_I^u=\sigma_{I_1}$ and $\sigma_J^u=\sigma_{J_1}$ where $I_1,J_1\in\Cal I$
satisfy one of the following conditions:
\roster
 \item "(a)" $I_1\cap I_1$ is a common endpoint of $I_1$ and $J_1$;
 \item "(b)" $\card(\partial I_1\cup\partial J_1)=4$.
\endroster
It is enough to exclude Case (b). Indeed, in this case
$\beta=\sigma_1^k\sigma_2^l$ cannot be conjugate to $\beta_1=\sigma_{I_1}^k\sigma_{J_1}^l$,
because $\lk_{i,j}(\beta^2)=0$ for $i\not\in\{1,2,3\}$ and any $j$ whereas
$\lk_{p,q}(\beta_1^2)=k$ and $\lk_{r,s}(\beta_1^2)=l$
for pairwise distinct $p,q,r,s$.
\qed\enddemo

\proclaim{ Lemma \lemGarside }
Let $X$ and $Y$ be two distinct conjugates of $\sigma_1$ in $\B_n$, $n\ge 3$.
If $XYX=YXY$, then there exists $u\in\B_n$ such that $X^u=\sigma_1$ and $Y^u=\sigma_2$.
\endproclaim

\demo{ Proof } Follows from Lemma \lemGarsideXYX.
\qed\enddemo

When speaking of Garside structures on groups, we use the terminology and notation from [\refQPtwo].
Let $(G,\Cal P,\delta)$ be a symmetric homogeneous square-free Garside structure
with set of atoms $\Cal A$ {\rm(for example, Birman-Ko-Lee's Garside structure [\refBKLone]
on the braid group, i.~e., $G=\B_n$, $\Cal A=\{\sigma_{ij}\}_{1\le i<j\le n}$,
$\Cal P=\{x_1\dots x_m\mid x_i\in\Cal A, m\ge 0\}$, $\delta=\sigma_{n-1}\sigma_{n-2}\dots\sigma_2\sigma_1$)}.

For $a,b\in G$ we set 
$b^G=\{b^a\mid a\in G\}$, and we write $a\sim b$ if $a\in b^G$ and
$a\preccurlyeq b$ if 
$a^{-1}b\in\Cal P$. We define the set of {\it simple elements} of $G$ as
$[1,\delta]=\{s\in G\mid 1\preccurlyeq s\preccurlyeq\delta\}$. For $X\in G$, the
{\it canonical length} of $X$ (denoted by $\ell(X)$) is the minimal $r$ such that
$X=\delta^p A_1\dots A_r$ for some $p\in\Z$, $A_1,\dots,A_r\in[1,\delta]\setminus\{1\}$.
The {\it summit length\/} of $X$ is defined as $\ell_s(X)=\min\{\ell(Y)\mid Y\in X^G\}$.
We denote the {\it cyclic sliding} of $X$ and the {\it set of sliding circuits} of $X$ by
$\frak s(X)$ and $\operatorname{SC}(X)$ respectively (these notions were introduced
in [\refGebGM], see also [\refQPtwo; Definition 1.12]).

\newpage
The following result was, in a sense, proven in [\refQPtwo] without stating it explicitly.

\proclaim{ Theorem \thQP }
Let $k,l\in\Z\setminus\{0\}$, $x,y\in\Cal A$, and let $Z=XY$
where $X\sim x^k$ and $Y\sim y^l$.
Then there exists $u\in G$ such that one of the following possibilities holds:
\roster
\item"(i)" $X^u=x_1^k$ and $Y^u=y_1^l$ with $x_1\in x^G\cap\Cal A$ and $y_1\in y^G\cap\Cal A$, or
\item"(ii)"
$\ell(Z^u)=\ell(X^u)+\ell(Y^u)$ and $Z^u\in\SC(Z)$.
\endroster
\endproclaim

\demo{ Proof } If the statement is true for $(k,l)$, then it is true for $(-k,-l)$,
therefore we may assume that $l>0$.
Then the proof of [\refQPtwo; Corollary 3.5] repeats almost word-by-word
in our setting if we define $\Cal Q_m$ as
$\{Z^u\mid u\in\Cal U_m\}$ where $\Cal U_m=\{u\mid\ell(X^u)\le 2m+|k|, \ell(Y^u)=l\}$.
Namely,
let $m$ be minimal under the assumption that $\Cal Q_m\ne\varnothing$. If $m=0$, then
(i) takes place. If $m>0$, then,
similarly to [\refQPtwo; Lemma 3.3] we show that if $u\in\Cal U_m$, then $\ell(Z^u)=\ell(X^u)+\ell(Y^u)$,
and
similarly to [\refQPtwo; Lemma 3.4] we show that $\frak s(\Cal Q_m)\subset\Cal Q_m$.
whence $\Cal Q_m\cap\SC(Z)\ne\varnothing$ which implies (ii).
\qed
\enddemo

\proclaim{ Corollary \corQP }
In the hypothesis of Theorem \thQP, assume that $Z$ is conjugate to $x^ky^l$.
Then there exists $u\in G$ such that (i) holds.
\endproclaim

\demo{ Proof } Suppose that (ii) takes place. Since
$\ell(Z)=\ell_s(Z)\le\ell(x^ky^l)\le |k|+|l|$, we have $\ell(X^u)+\ell(Y^u)\le |k|+|l|$.
By combining this fact with $\ell(X^u)\ge |k|$ and $\ell(Y^u)\ge |l|$, we obtain
$\ell(X^u)=|k|$ and $\ell(Y^u)=|l|$, and the result follows from [\refQPtwo; Theorem 1a].
\qed\enddemo

\proclaim{ Lemma \lemGarsideXYX }
Let $X\sim x$ and $Y\sim y$ where $x,y\in\Cal A$.
If $XYX=YXY$, then there exists $u\in G$ such that $X^u,Y^u\in\Cal A$.
\endproclaim

\demo{ Proof }
Without loss of generality we may assume that $Y=y\in\Cal A$.
By [\refQPtwo; Theorem 1a], the left normal form of $X$ is
$\delta^{-p}\cdot A_p\cdot\dots\cdot A_1\cdot x\cdot B_1\cdot\dots\cdot B_p$
where $A_i,B_i\in[1,\delta]\setminus\{1\}$, $A_i\delta^{i-1}B_i = \delta^i$ for $i=1,\dots,p$.
By symmetry, the right normal form of $X$ is
$C_p\cdot\dots\cdot C_1\cdot x\cdot D_1\cdot\dots\cdot D_p\cdot\delta^{-p}$
again with  $C_i,D_i\in[1,\delta]\setminus\{1\}$, $C_i\delta^{i-1}D_i = \delta^i$ for $i=1,\dots,p$.
We have $\sup yXy\le 2\sup y + \sup X = 3+p$.
Thus, if $p>0$, then $\sup X + \sup y + \sup X = 3+2p > 3+p = \sup xYx = \sup YxY$.
Then, by [\refQPii; Lemma 2.1b], either $B_p y$ or $yC_p$ is a simple element.
Without loss of generality we may assume that $B_p y\in[1,\delta]$.

Since the Garside structure is symmetric and $B_p y$ is simple,
there exists an atom $y_1$ such that
$B_p y=y_1 B_p$. Thus, for $v=B_p^{-1}$, we have $y^v\in\Cal A$ and the left normal form
of $X^v$ is $\delta^{p-1}\cdot A_{p-1}\cdot\dots\cdot A_1\cdot x\cdot B_1\cdot\dots\cdot B_{p-1}$.
Therefore, the induction on $p$ yields $X^u=x$ and $Y^u=z\in\Cal A$ for $u=(B_1\dots B_p)^{-1}$.
\qed\enddemo

\medskip\noindent
{\bf Remark \remGarside. } (Compare with  [\refCP, \refLM]). Lemma \lemGarsideXYX\
admits the following ge\-ne\-ra\-li\-za\-tion which can be proven using
the results and methods of [\refQPtwo, \refQPii].
{\sl Let $X\sim x$ and $Y\sim y$ for $x,y\in\Cal A$.
Then either $X$ and $Y$ generate a free subgroup of $G$, or
there exists $u\in G$ such that $X^u,Y^u\in\Cal A$.
In the latter case, the subgroup generated by $X,Y$ is either free or
isomorphic to Artin group of type $I_2(p)$, $p\ge 2$.}
In particular, for $G=\B_n$, {\sl if $X$ and $Y$ are two conjugates of $\sigma_1$,
then either $X$ and $Y$ generate a free subgroup of $\B_n$, or
there exists $u\in\B_n$ such that $X^u=\sigma_1$ and $Y^u=\sigma_i$ for some $i$.}
Maybe, I will write a proof of this fact in a future paper. 

\subhead Acknowledgement \endsubhead
I am grateful to the referee for many very useful remarks.

\Refs
\def\r{\ref}

\r\no\refArtin
\by E.~Artin \paper Theory of braids \jour Ann. of Math. \vol 48 \yr 1947 \pages 101--126
\endref

\r\no\refBLM
\by J.S.~Birman, A.~Lubotzky, J.~McCarthy
\paper Abelian and solvable subgroups of the mapping class group
\jour Duke Math. J. \vol  50 \yr 1983 \pages 1107--1120
\endref

\r\no\refBKLone
\by J.~Birman, K.-H.~Ko, S.-J.~Lee
\paper A new approach to the word and conjugacy problems in the braid groups
\jour  Adv. Math. \vol 139 \yr 1998 \pages 322--353
\endref

\r\no\refCP
\by    J.~Crisp, L.~Paris
\paper The solution to a conjecture of Tits on the subgroup generated by the squares of
       the generators of an Artin group \jour Invent. Math. \vol 145 \yr 2001 \pages 19--36
\endref

\r\no\refDG
\by    J.L.~Dyer, E.K.~Grossman
\paper The automorphism group of the braid groups
\jour  Amer. J. of Math. \vol 103 \yr 1981 \pages 1151--1169
\endref

\r\no\refEil
\by    S. Eilenberg
\paper Sur les transformations p\'eriodiques de la surface de sph\`ere,
\jour  Fund. Math. \vol 22 \yr 1934 \pages 28--41
\endref

\r\no\refFH
\by    W.~Fulton, J.~Harris
\book  Representation theory. A first course
\publ  Springer \yr 1991
\endref

\r\no\refGebGM
\by V.~Gebhardt, J.~Gonz\'alez-Meneses
\paper The cyclic sliding operation in Garside groups
\jour  Math. Z. \vol 265 \yr 2010 \pages 85--114
\endref

\r\no\refGM
\by    J.~Gonz\'alez-Meneses
\paper The $n$th root of a braid is unique up conjugacy
\jour  Algebraic and Geometric Topology \vol 3 \yr 2003 \pages 1103--1118
\endref

\r\no\refGMW
\by    J.~Gonz\'alez-Meneses, B.~Wiest
\paper On the structure of the centralizer of a braid
\jour  Ann. Sci. \'Ec. Norm. Sup\'er. (4) \vol 37 \yr 2004 \pages 729--757
\endref

\r\no\refGL
\by    E.A.~Gorin, V.Ya.~Lin,
\paper Algebraic equations with continuous coefficients and some problems
       of the algebraic theory of braids
\jour  Math. USSR-Sbornik \vol 7 \yr 1969 \pages 569–596.
\endref

\r\no\refI
\by    N.~V. Ivanov
\book  Subgroups of Teichm\"uller modular groups
\bookinfo Translations of mathematical monographs \vol 115 \yr 1992 \publ AMS
\endref

\r\no\refKer
\by    B. de Ker\'ekj\'art\'o
\paper \"Uber die periodischen Transformationen der Kreisscheibe und
       der Kugel\-fl\"a\-che \jour Math. Annalen \vol 80 \yr 1919 \pages 3--7
\endref

\r\no\refLM
\by    C.~Leininger, D.~Margalit
\paper Two-generator subgroups of the pure braid group
\jour  Geom. Dedicata \vol 147 \yr 2010 \pages 107--113
\endref

\r\no\refLin
\by    V.~Lin \paper Braids and permutations \jour arXiv:math/0404528 \endref

\r\no\refLinTalk
\by    V.~Lin
\paper Some problems that I would like to see solved
\jour  Abstract of a talk. Technion, 2015,
http://www2.math.technion.ac.il/$\widetilde{\;}$pincho/Lin/Abstracts.pdf
\endref

\r\no\refMMP
\by    K.~Magaard, G.~Malle, P.H.~Tiep
\paper Irreducibility of tensor squares, symmetric squares and alternating squares
\jour  Pac. J. Math. \vol 202 \yr 2002 \pages 379-427
\endref

\r\no\refMKS
\by    W.~Magnus, A.~Karrass, D.~Solitar
\book  Combinatorial group theory: presentations of groups in terms of generators and relations
\publ  Interscience Publ. \yr 1966
\endref

\r\no\refM
\by    S.~Manfredini
\paper Some subgroups of Artin’s braid group. Special issue on braid
groups and related topics (Jerusalem, 1995) \jour Topology Appl. \vol 78 \yr 1997 \pages 123--142
\endref


\r\no\refUR
\by     S.Yu.~Orevkov
\paper  Quasipositivity test via unitary representations of braid groups
        and its applications to real algebraic curves
\jour   J. Knot Theory Ramifications \vol 10 \yr 2001  \pages 1005--1023
\endref

\r\no\refQPtwo
\by     S.Yu.~Orevkov
\paper  Algorithmic recognition of quasipositive braids of algebraic length two
\jour   J. of Algebra \vol 423 \yr 2015 \pages 1080--1108
\endref

\r\no\refQPii
\by     S.Yu.~Orevkov
\paper  Algorithmic recognition of quasipositive 4-braids of algebraic length three
\jour   Groups, Complexity, Cryptology \vol 7 \yr 2015 \issue 2 \pages 157--173 
\endref

\endRefs
\enddocument